\documentclass[12pt]{amsart}
\usepackage{latexsym, amssymb, amsmath, amscd}
 \usepackage{eucal}

    \newtheorem{rema}{Remark}[section]
    \newtheorem{propo}[rema]{Proposition}

   \newtheorem{theo}[rema]{Theorem}
   \newtheorem{def-theo}[rema]{Definition-Theorem}
 
   \newtheorem{defi}[rema]{Definition}
    \newtheorem{lemma}[rema]{Lemma}
    \newtheorem{corol}[rema]{Corollary}
     
  \newtheorem{rmk}[rema]{Remark}

    \newtheorem{prob}[rema]{Problem}

	\newcommand{\nno}{\nonumber}

 \newcommand{\pf}{{\it Proof:}\hspace{2ex}}
 
 \newcommand{\epfv}{\hspace{1em}$\Box$\vspace{1em}}

\newcommand{\bC}{{\mathbb C}}
\newcommand{\bZ}{{\mathbb Z}}

\newcommand{\cA}{{\mathcal A}}

\newcommand{\cM}{{\mathcal M}}
\newcommand{\cN}{{\mathcal N}}
\newcommand{\cB}{{\mathcal B}}

\newcommand{\vt}{\vartheta} 
\newcommand{\vms}{$\vartheta$-Mathieu subspace }
\newcommand{\vmss}{$\vartheta$-Mathieu subspaces }

\title[Mathieu Subspaces of Modules of Associative Algebras]
{A Generalization of Mathieu Subspaces to Modules of Associative Algebras} 

  \author{Wenhua Zhao}      
\address{ Department of Mathematics, Illinois State University, Normal, IL 61790-4520. 
{\em E-mail}: wzhao@ilstu.edu.}

\begin{document}

\begin{abstract}
We first propose a generalization of the notion of Mathieu subspaces of associative algebras $\cA$, which was introduced recently in \cite{GIC} and \cite{MS}, to $\cA$-modules $\cM$. The newly introduced notion in a certain sense also generalizes the notion of submodules. Related with this new notion, we also introduce 
the sets $\sigma(N)$ and $\tau(N)$ of {\it stable elements} and 
{\it quasi-stable elements}, respectively, for all $R$-subspaces 
$N$ of $\cA$-modules $\cM$, where $R$ is the base ring of $\cA$. 
We then prove some general properties of the sets 
$\sigma(N)$ and $\tau(N)$. Furthermore, examples from certain modules of the {\it quasi-stable algebras} \cite{MS}, matrix algebras over fields and polynomial algebras 
are also studied. 
\end{abstract}

\keywords{Mathieu subspaces of associative algebras, Mathieu subspaces of modules of associative algebras, (quasi-)stable elements, (quasi-)stable algebras, (quasi-)stable modules.}

\subjclass[2000]{16D10, 16D99}

\thanks{The author has been partially supported 
by NSA Grant H98230-10-1-0168}

 \bibliographystyle{alpha}
    \maketitle


\renewcommand{\theequation}{\thesection.\arabic{equation}}
\renewcommand{\therema}{\thesection.\arabic{rema}}
\setcounter{equation}{0}
\setcounter{rema}{0}
\setcounter{section}{0}

\section{\bf Introduction}\label{S1}

\subsection{Background and Motivation}
\label{S1.1}
Once and for all, we fix an arbitrary unital commutative ring 
$R$ and an associative 
but not necessarily commutative unital 
$R$-algebra $\cA$. Recall that 
the following notion has been introduced 
recently in \cite{GIC} and \cite{MS}. 
 
\begin{defi}\label{Def-MS}
Let $J$ be a  $R$-submodule or 
$R$-subspace of $\cA$. 
We say $J$ is a {\it left} 
$($resp.,  {\it right}; {\it two-sided}$)$ 
{\it  Mathieu subspace} of $\cA$  
if the following property holds: 
if $a\in \cA$ satisfies $a^m\in J$ for all $m\ge 1$,
then for each $b, c\in \cA$, we have 
$ b a^m \in J$ $($resp., $a^m c \in J$; $ba^m c \in J$$)$ 
for all $m\gg 0$, i.e., there exists $N\ge 1$ 
such that $b a^m  \in J$ $($resp.,  $a^mc \in J$; $ba^m c \in J$$)$ 
for all $m\ge N$.
\end{defi}   

Two-sided Mathieu subspaces will also simply be 
called Mathieu subspaces. A $R$-subspace $M$ of $\cA$ is said to be a {\it pre-two-sided} Mathieu subspace of $\cA$ if it is both left and right Mathieu subspace of $\cA$. Note that the {\it pre-two-sided} Mathieu subspaces were previously called  {\it two-sided} Mathieu subspace 
or {\it Mathieu subspaces} in \cite{GIC}.

The introduction of the notion of Mathieu subspaces in \cite{GIC} 
and \cite{MS} was mainly motivated by the studies of 
{\it the Jacobian conjecture} \cite{K} 
(see also \cite{BCW} and \cite{E}), 
{\it the Dixmier conjecture} \cite{D} 
(see also \cite{Ts}, \cite{BK} and \cite{AE}), 
{\it the Mathieu conjecture}  \cite{Ma}, 
{\it the vanishing conjecture} \cite{HNP}, \cite{GVC}, 
\cite{AGVC}, \cite{EWiZ} and more recently, {\it the image conjecture} \cite{IC},  
as well as many other related open problems. 
Actually, both {\it the Mathieu conjecture} and {\it the image conjecture} imply 
{\it the Jacobian conjecture} (and hence also {\it the Dixmier conjecture}), 
and both are (open) problems on whether or not certain subspaces 
of some algebras are Mathieu subspaces.  
For some recent developments on Mathieu subspaces, 
see \cite{MS}, \cite{FPYZ}, \cite{EWZ1}, \cite{EWZ2}, \cite{WZ} 
and \cite{EZ}.    
For a recent survey on the {\it the image conjecture} 
and it's connections with some other problems, 
see \cite{E2}.  

Note that every left ideal of $\cA$ is a left 
Mathieu subspace of $\cA$. This is also the case for right or two-sided ideals. Therefore, the notion of Mathieu subspaces can be viewed as a natural generalization of the notion of  ideals. 

On the other hand, the notion of ideals of  associative algebras $\cA$ has a natural generalization to their modules.   
Namely, viewing $\cA$ itself as a {\it left}  
(resp., right or two-sided) $\cA$-module in the canonical way, every {\it left} 
(resp., {\it right}  or {\it two-sided}) 
ideal of $\cA$ is a submodule of the 
{\it left} (resp., {\it right} or {\it two-sided}) 
$\cA$-module $\cA$. 

Naturally, one may wonder whether or not there is also a generalization of Mathieu subspaces of 
associative algebras to their  
modules, which is parallel 
to the above generalization of the 
ideals to submodules, such that the 
following diagram of notions commutes: 
\allowdisplaybreaks{
\begin{align}\label{NotionCD} 
\begin{CD}
{\bf  Ideal}\, @> \mathcal G_1>> \quad {\bf Submodule \quad} \\
@V \mathcal G_2 VV @     VV\mathcal G_4 V \\
\mbox{\bf \quad Mathieu Subspace}
\quad @ >\mathcal G_3 >> \, {\bf ??}
\end{CD}
\end{align} 
where the ``maps" $\mathcal G_i$ $(1\le i\le 4)$  in the diagram denote the generalizations 
of the corresponding notions.  }

Note that the generalization from the notion of  ideals to the notion of submodules is crucial in the theories of various algebras (associative algebras, Lie algebras, etc.). We believe that the generalization $\mathcal G_3$ in the diagram above is also important, not only for the study of modules of associative algebras but also for the study of Mathieu subspaces of associative algebras. 

In this paper, we first complete the 
commutative diagram above by introducing 
the notion of what we call {\it Mathieu subspaces} 
for modules $\cM$ of associative algebras 
$\cA$ (see Definition \ref{Def-GMS} below).
Related with this new notion, we also introduce 
the sets $\sigma(N)$ and $\tau(N)$ of {\it stable}  elements and {\it quasi-stable} elements, respectively, for all $R$-subspaces $N$ of $\cA$-modules $\cM$ (see the discussions in Subsection \ref{S1.2} below and also 
Definition \ref{Def-Quasi-Stable} 
in Section \ref{S3}).  We then study some general properties of the sets $\sigma(N)$ and $\tau(N)$. Furthermore,  examples from certain modules of the {\it quasi-stable algebras} \cite{MS} (see also Definition \ref{q-StaAlg}), matrix  algebras over fields and polynomial algebras are also studied. 

\subsection{To Complete the Commutative Diagram} 
\label{S1.2} 
To be more precise, let's first recall the following 
terminologies and conventions introduced in \cite{MS}, 
which will also be used throughout this paper.  

Note first that for Mathieu subspaces and ideals, 
as well as for the new notions to be defined below,  
we have several different cases:  
{\it left, right} and {\it $($pre-$)$two-sided}. 
Very often, it is necessary and important to treat all 
these cases. For simplicity, following \cite{MS} we introduce 
the short terminology {\it $\vartheta$-Mathieu subspaces} for 
Mathieu subspaces, where $\vartheta$ stands for 
{\it left, right}, {\it pre-two-sided}, or 
{\it two-sided}. Similarly, we introduce 
the terminology {\it $\vartheta$-ideals} for ideals, 
except for the specification  
$\vartheta=${\it ``pre-two-sided\hspace{.1mm}"}, 
we also set {\it $\vartheta$-ideals} to mean 
{\it two-sided} ideals.    

In other words, the reader should read the letter 
$\vartheta$ as an index or a variable with 
four possible choices or 
``values".  However, to avoid repeating the phrase 
``{\it for every specification of $\vartheta$}" 
or {\it ``for every $\vt$"} 
infinitely many times, we will simply 
leave $\vartheta$ unspecified 
for the statements or propositions which hold 
for all the four specifications of $\vartheta$.  
For example, by saying ``{\it Any  
$\vt$-ideal of $\cA$ is a \vms of $\cA$}" 
or ``{\it $\cA$ has at least one proper 
$\vartheta$-Mathieu subspace}",  
we mean that the statement holds for every 
fixed specification of $\vartheta$: 
{\it left, right}, {\it pre-two-sided} 
or {\it two-sided}.  

Furthermore, in the case that only the 
{\it two-sided} case is concerned, 
the variable $\vt$ will be simply 
dropped, e.g.,  the words {\it ideals} and 
{\it Mathieu subspaces} (without 
the variable $\vt$) always mean 
the  two-sided ideals and two-sided 
Mathieu subspaces, respectively.

%
%

Now, let $R$ and $\cA$ be as before and $\cM$ 
a left $\cA$-module. For any $u\in \cM$ and 
any subset $N$ of $\cM$, we set 
\begin{align}\label{Def-Ann}
(N:u)\!:=\{a\in \cA\,|\, au\in N\}.
\end{align}

With the terminologies and notations fixed above, 
we can now introduce the new notion claimed 
in the previous subsection as follows. 

\begin{defi}\label{Def-GMS}
Let $u\in \cM$ and $N$ be a $R$-subspace or $R$-submodule 
of $\cM$. We say that $N$ is a \vms of $\cM$ 
with respect to $u$ if $(N:u)$ 
is a \vms of $\cA$.  
\end{defi}

In a similar way as above we can define \vmss for right $\cA$-modules as well. 
Since every right $\cA$-module 
is a left $\cA^{op}$-module, 
where $\cA^{op}$ denotes {\it the opposite algebra} 
of $\cA$, 
we may (and also will) focus only on 
the \vmss of left $\cA$-modules. 
Therefore, throughout this paper, 
all $\cA$-modules will be assumed 
to be left $\cA$-modules unless 
stated otherwise. In particular, 
when we view $\cA$ itself as an 
$\cA$-module, we always mean 
the left $\cA$-module over 
the $R$-space $\cA$ with 
the left action given by 
the algebra product of $\cA$.
  
\vspace{2mm}
The main ideas behind Definition 
\ref{Def-GMS} are as follows. 

\vspace{2mm}
Note first that in general 
$\cA$-modules $\cM$ have no product operations. Therefore, 
we cannot generalize directly the notion of \vmss of the algebra 
$\cA$ to $\cM$, as we did for the generalization 
of left ideals of $\cA$ to submodules of $\cM$.   

But, on the other hand, for every $R$-subspace 
$J$ of the $\cA$-module $\cA$, 
we have $(J:1_\cA)=J$. Hence, $J$ is a \vms of the algebra $\cA$ iff $J$ is a \vms of the 
$\cA$-module $\cA$ with respect to 
$1_\cA\in \cA$. Therefore, from this point of view we see that the notion of {\it Mathieu subspaces} of $\cA$-modules given in Definition \ref{Def-GMS}
indeed generalizes the notion of {\it Mathieu subspaces} of the algebra $\cA$. 
Actually, as we will show later in Section \ref{S3} 
(see Corollary \ref{SubmoduleLemma} and 
Remark \ref{rmk-cG4}), with $\vt=${\it ``left"} it also generalizes 
the notion of submodules of $\cA$-modules $\cM$ in the sense that: 
{\it for every $R$-subspace $N$ of 
$\cM$, $N$ is a submodule of $\cM$ 
iff $N$ is a left Mathieu subspace of $\cM$ 
with respect to all elements $u\in \cM$}.     
 
In contrast to ideals of $\cA$ or submodules of $\cA$-modules 
$\cM$, the \vmss $N$ of $\cA$-modules $\cM$ defined in Definition \ref{Def-GMS} depend on some {\it referring} elements $u\in \cM$. So it is important and also convenient to consider the set of all the {\it referring} elements of 
$N$, i.e., the set of all the elements $u\in \cM$ with respect to which $N$ is a \vms  of $\cM$. We denote this set by $\tau_\vt(N)$  \label{TWO-SETS} 
and call it the set of {\it $\vt$-quasi-stable} elements 
of the $R$-subspace $N\subseteq \cM$. 
Similarly, we let $\sigma_\vt(N)$ denote the set  
of all the elements $u\in \cM$ such that 
$(N:u)$ is a $\vt$-ideal of $\cA$, and call it  
the set of {\it $\vt$-stable} elements 
of $N$. For the justifications of the 
terminologies of {\it $\vt$-quasi-stableness} and 
{\it $\vt$-stableness}, see Definition 
\ref{Def-Quasi-Stable} and also the followed  
discussions.


Even though our main concerns in this paper are 
on the sets $\tau_\vt(N)$ of $\vt$-quasi-stable 
elements of $R$-subspaces $N$ of $\cA$-modules $\cM$, 
for the purpose of comparison we also give 
all the parallel results (to that of 
$\tau_\vt(N)$) for the set $\sigma_\vt(N)$ 
of $\vt$-stable elements of $N$. 
But, since the proofs for the case of 
$\sigma_\vt(N)$ are very often similar to  
(and actually simpler than) the proofs for 
the case of $\tau_\vt(N)$, we will leave 
to the reader most of the proofs for 
the case of $\sigma_\vt(N)$.

\subsection{Arrangement} \label{S1.3}
In Section \ref{S2}, we recall some results on 
\vmss of associative algebras $\cA$, which will be 
needed later in this paper. Some of these results 
will also be generalized in Section \ref{S3} to 
\vmss of $\cA$-modules.  

In Section \ref{S3}, we study some general 
properties of the sets $\tau_\vt(N)$ and $\sigma_\vt(N)$ 
of {\it $\vt$-quasi-stable} elements and {\it $\vt$-stable}  
elements, respectively, of $R$-subspaces $N$ of 
$\cA$-modules $\cM$. 

Note that the processes of taking the sets  
$\tau_\vt(N)$ and $\sigma_\vt(N)$  
of $R$-subspaces $N$ of $\cM$ 
can be viewed as two maps from the set  
of $R$-subspaces $N$ of $\cM$ to 
the set of subsets of $\cM$. 
In this section, we mainly study the relations of 
these two maps with some other operations 
on $R$-subspaces of $\cM$, 
such as the relation with the operation 
of finite intersections in Proposition \ref{wIntersectn}; 
the relation with the pulling-back by the homomorphisms 
of $\cA$-modules in Proposition \ref{Inverse-vee-1}; 
and the relation with the pulling-back by homomorphisms of algebras in Proposition \ref{Inverse-vee-2}, etc. 
Some consequences of the propositions above 
are also derived.

Another main result of this section is Theorem 
\ref{Max-submodule}, which claims that for every 
$R$-subspace $N$ of $\cM$, the intersections of both 
$\tau_\vt(N)$ and $\sigma_\vt(N)$ with $N$ itself 
are the same as the (unique) $\cA$-submodule 
of $\cM$, denoted by $I_N$, which is maximum among 
all the $\cA$-submodules 
of $\cM$ contained in $N$. For short,  
in this paper we call the $\cA$-submodule $I_N$ 
the {\it maximum $\cA$-submodule of $N$}.    
Some consequences of this theorem are also 
derived in Corollaries \ref{Corol4Irred}--\ref{SubmoduleLemma}.

In Section \ref{S4}, we first introduce 
in Definition \ref{Def-QSM} the notions of  
the {\it $\vt$-quasi-stable} and {\it $\vt$-stable}  
$\cA$-modules. The relations of these two new notions 
with the notions of {\it $\vt$-quasi-stable} and 
{\it $\vt$-stable} $R$-algebras introduced in \cite{MS} 
(see also  Definition \ref{q-StaAlg}) are given 
in Proposition \ref{qStableEquiv}. More explicitly,     
a  $R$-algebra $\cA$ is {\it $\vt$-quasi-stable} 
(resp., {\it $\vt$-stable}) iff all its modules are 
{\it $\vt$-quasi-stable} (resp., {\it $\vt$-stable}). 

One special property of {\it $\vt$-quasi-stable} 
(resp., {\it $\vt$-stable}) modules $\cM$ is that 
for all $R$-subspaces $N$ of $\cM$, 
by Theorem \ref{Max-submodule} we have $\tau_\vt(N)=I_N\cup N^c$ 
(resp., $\sigma_\vt(N)=\tau_\vt(N)=I_N\cup N^c$), 
where $I_N$ is the {\it maximum $\cA$-submodule of $N$} 
(defined above) and $N^c$ denotes the {\it complement} of $N$ in $\cM$. 

On the other hand, it has been shown in \cite{MS} 
(see also Proposition \ref{q-Quasi/R}) that for every 
$R$-algebra $\cA$, if $\cA$ is integral over $R$ and every  
element of $\cA$ is either invertible or nilpotent, 
then $\cA$ is quasi-stable. 
Consequently, every left or right Artinian $R$-algebra $\cA$, 
which is integral over $R$, is quasi-stable. 
The {\it $\vt$-quasi-stable} and 
{\it $\vt$-stable} algebras over fields have been 
classified in \cite{MS} (see also Theorem 
\ref{Class-qStable} and Propositions \ref{AlgLocal} 
and \ref{Class-Stable}). 

Therefore, we get a family of associative algebras 
$\cA$ such that for all $\cA$-modules $\cM$, the sets 
$\tau_\vt(N)$ for all $R$-subspaces $N$ of $\cM$ 
can be determined explicitly 
(see Corollaries \ref{Tau4qStable} 
and \ref{DivisionCase}). 

In Section \ref{S5}, we study the cases for two modules of    
the matrix algebras $M_n(K)$ $(n\ge 1)$ 
over an arbitrary field $K$. First, we determine in Proposition \ref{KnCase} the sets $\tau_\vt(N)$ and $\sigma_\vt(N)$ explicitly for all $K$-subspaces $N$ of the $M_n(K)$-module  
$K^n$ (with the canonical left action). Second, we view 
$M_n(K)$ itself as a (left) $M_n(K)$-module 
in the canonical way and determine 
in Proposition \ref{V-CoD1} the sets $\tau_\vt(N)$ 
and $\sigma_\vt(N)$ explicitly for all co-dimension one 
$K$-subspaces $N$ of $M_n(K)$.

In Section \ref{S6}, motivated by Conjecture $3.2$ 
in \cite{GIC} on integrals of polynomials, 
we study the cases for two modules of the polynomial algebras 
over fields $K$. We first consider a family of $K$-subspaces 
$N_{B, \alpha}$ defined in Eq.\,(\ref{Def-NB}) 
of the polynomial algebra in several variables 
over an arbitrary field $K$, and derive the sets 
$\sigma(N_{B, \alpha})$ and $\tau(N_{B, \alpha})$  
explicitly in Proposition \ref{Nf-propo}. 
We then consider a family of $\bC$-subspaces 
$N_q$ defined in Eq.\,(\ref{Def-Nq}) of the polynomial algebra in one variable over the complex field $\bC$, 
and derive the sets $\sigma(N_q)$ and $\tau(N_q)$   explicitly in Proposition \ref{NqCase}.

\renewcommand{\theequation}{\thesection.\arabic{equation}}
\renewcommand{\therema}{\thesection.\arabic{rema}}
\setcounter{equation}{0}
\setcounter{rema}{0}

\section{\bf Some Properties of $\vt$-Mathieu Subspaces of Associative Algebras} \label{S2}

Let $R$ and $\cA$ be as fixed in the previous 
section. We denote by $1_\cA$ or simply $1$ the 
identity element of the $R$-algebra $\cA$.  
The sets of units or 
invertible elements of $R$ and 
$\cA$ will be denoted by $R^\times$ and 
$\cA^\times$, respectively.  
All the notations, conventions and terminologies fixed 
in the previous section will also be in force throughout 
the paper.

Let $J$ be an arbitrary subset of $\cA$. 
Following \cite{MS} we define the {\it radical} of $J$, 
denoted by $\sqrt{J}$, 
to be $\sqrt{J}\!:=\{a\in \cA\,|\, a^m\in J \mbox{ when } m\gg 0\}$. 
Note that $\sqrt{\cA}=\cA$, and $\sqrt 0$ is the set of 
all nilpotent elements of $\cA$. We will also denote $\sqrt 0$  
by ${\rm nil\,}(\cA)$ since when $\cA$ is commutative, 
$\sqrt 0$ is the same as the {\it nilradical} of $\cA$. 

For any $a\in \cA$ and $\vt\neq${\it ``pre-two-sided"}, 
we denote by $(a)_\vt$ the $\vt$-ideal generated by $a$. 
For the case $\vt=${\it ``pre-two-sided"}, we set $(a)_\vt\!:=Aa+aA$. 
Also, as fixed in the previous section, when 
$\vt=${\it ``two-sided"}, we will simply drop the 
variable $\vt$ from the notation above, 
i.e., $(a)$ denotes the (two-sided) 
ideal of $\cA$ generated by the element $a$.

Furthermore, for each $a\in \cA$, we also set 
\begin{align} 
a^{-1}J\!:&=\{b\in \cA\,|\, ab \in J\}.     
\label{a-Inv-J}
\end{align}
Note that $a^{-1}J$ is an abusing notation 
since $a$ might not be invertible in $\cA$. 

In this section, we mainly recall some properties of \vmss of $R$-algebras $\cA$, which will be needed later in this paper.

Let's start with the following lemma, which is very simple but provides 
a family of (two-sided) Mathieu subspaces.

\begin{lemma}\label{radical-Lemma}
Let $J$ be a  $R$-subspace of $\cA$ such that 
$\sqrt{J}\subseteq {\rm nil\,}(\cA)$. Then 
every $R$-subspace $H$ of $J$ is a 
Mathieu subspace of $\cA$. 
\end{lemma} 

The next lemma proved in \cite{MS} provides 
a different point of view to see that 
the notion of \vmss is indeed a natural 
generalization of that of $\vt$-ideals.

\begin{lemma}\label{CharsInRads}
Let $J$ be a  $R$-subspace of a  $R$-algebra 
$\cA$. Then the following statements hold.
 \begin{enumerate}
   \item[$i)$] $J$ is a left ideal of $\cA$ iff 
 for any $a\in \cA$, we have 
  \begin{align*}
J \subseteq a^{-1}J.
\end{align*}

  \item[$ii)$] $J$ is a left Mathieu subspace of $\cA$ iff 
 for any $a\in \cA$, we have 
  \begin{align*} 
\sqrt{J} \subseteq  \sqrt{a^{-1}J}.
\end{align*}

  \item[$iii)$] $J$ is a right ideal of $\cA$ iff 
for any $a\in \cA$, we have
  \begin{align*}
J\subseteq (J: a).
\end{align*}

  \item[$iv)$] $J$ is a right Mathieu subspace of 
$\cA$ iff for any $a\in \cA$, we have 
\begin{align*} 
\sqrt{J} \subseteq \sqrt{(J: a)}.
\end{align*}

  \item[$v)$] $J$ is a $($two-sided$)$ ideal of $\cA$ iff 
for any $a, b\in \cA$, we have
  \begin{align}\label{CharsInRads-e7}
J\subseteq (a^{-1}J: b).
\end{align}

\item[$vi)$] $J$ is a $($two-sided$)$ Mathieu subspace of 
$\cA$ iff for any $a, b\in \cA$, we have 
  \begin{align}\label{CharsInRads-e8}
\sqrt{J} \subseteq \sqrt{(a^{-1}J: b)}.
\end{align}
\end{enumerate} 
\end{lemma}

The following two simple lemmas, first noticed in \cite{GIC}, 
will be very useful for our later arguments.

\begin{lemma}\label{ONE-Lemma}
Any proper $R$-subspace $J$ of $\cA$ 
with $1_\cA\in J$ is not 
a \vms of $\cA$.   
\end{lemma}

\begin{lemma}\label{IntersectionLemma}
For every finite family of \vmss 
$J_i$ $(1\le i\le k)$ of $\cA$, the intersection 
$\bigcap_{i=1}^k J_i$ is also a \vms of $\cA$.  
\end{lemma}

Note that in contrast to the case of $\vt$-ideals, 
intersections of infinitely many \vmss of $\cA$ are 
not always \vmss of $\cA$, e.g., see Example 
$4.17$ in \cite{MS}. 

The following two propositions  will also be 
important for our later arguments. 

\begin{propo}\label{J-quotient} 
$($\cite{MS}$)$\,
Let $I$ be an ideal of $\cA$ and 
$J$ a  $R$-subspace of $\cA$ 
which contains $I$. Then $J$ is a 
\vms of $\cA$ iff $J/I$ 
is a \vms of $\cA/I$.
\end{propo}

\begin{propo}\label{Pull-Back} $($\cite{GIC}$)$\,  
Let $\phi:\cA\to \cB$ be a homomorphism of 
$R$-algebras. Then for every \vms $N$ of $\cB$, 
$\phi^{-1}(N)$ is a \vms 
of $\cA$.
\end{propo} 

%

Next, we recall the following notions introduced in 
\cite{MS}.  

\begin{defi}\label{q-StaAlg}
A $R$-algebra $\cA$ is said to be 
{\it $\vartheta$-quasi-stable} $($resp., {\it $\vartheta$-stable}$)$ 
if every $R$-subspace $J$ of $\cA$ with $1\not \in J$ is a 
\vms $($resp., $\vartheta$-ideal$)$ of $\cA$.
\end{defi}

A family of quasi-stable $R$-algebras is given by the following 
proposition which was proved in Proposition $7.4$ and 
Corollary $7.5$ in \cite{MS}. 

\begin{propo}\label{q-Quasi/R}
Let $\cA$ be a $R$-algebra such that $\cA$ is integral over $R$ 
and every element of $\cA$ is either invertible or 
nilpotent. Then $\cA$ is a quasi-stable $R$-algebra.  

Consequently, every left or right Artinian 
local $R$-algebra, which is integral over $R$,  
is quasi-stable. 
\end{propo}

Now, we assume that the base ring $R$ is a field 
$K$ and $\cA$ is a  $K$-algebra. 
A subset $S\subseteq \cA$ is said to be {\it algebraic} 
over $K$ if every element $a\in S$ 
is algebraic over $K$, i.e., $a$ satisfies 
a nonzero polynomial in one variable with 
coefficients in $K$. Recall also that 
an element $a\in \cA$ is said to 
be an {\it idempotent} if $a^2=a$.   

With the terminologies fixed above, we have the following characterization 
proved in \cite{MS} for the \vmss $J$ of $K$-algebras $\cA$, 
whose radical $\sqrt J$ is algebraic over $K$. 

\begin{theo} \label{CharByIdem}
Let $J$ be a  $K$-subspace of a  $K$-algebra $\cA$ 
such that $\sqrt J$ is algebraic over $K$. 
Then $J$ is a \vms of $\cA$ iff 
for every idempotent $e\in J$, we have 
$(e)_\vt\subseteq J$. 
\end{theo}

%


Furthermore, the {\it $\vartheta$-quasi-stable} 
algebras over fields $K$ defined above 
in Definition \ref{q-StaAlg} 
have been classified in \cite{MS} 
as follows.

\begin{theo}\label{Class-qStable}
Let $K$ be an arbitrary field and $\cA$ a  $K$-algebra. 
Then $\cA$ is $\vartheta$-quasi-stable 
iff either $\cA \simeq K \dot{+} K$, or $\cA$ is 
an algebraic local $K$-algebra, where 
$K \dot{+} K$ is the $K$-algebra with 
$K\times K$ as the base $K$-vector space 
and the component-wise product as 
the algebra product.   
\end{theo}

One remark on the theorem above is that 
for algebras over fields, 
the $\vartheta$-quasi-stableness actually does 
not depend on the specialization of $\vt$. 
More precisely, an algebra $\cA$ over a field $K$ 
is one-sided (left or right) quasi-stable iff 
it is (two-sided) quasi-stable. 
By Proposition \ref{Class-Stable} below, 
this is also the case for the 
$\vartheta$-stableness.
 
\begin{propo}\label{AlgLocal} 
For any algebraic $K$-algebra $\cA$, the following statements 
are equivalent:

$1)$ $\cA$ is local; 

$2)$ $\cA$ has no idempotents except $0, 1 \in \cA$;

$3)$ each element of $\cA$ is either nilpotent 
or invertible. 
\end{propo} 

Finally, let's point out that the $\vt$-stable $K$-algebras  
have also been classified in \cite{MS} as follows.

\begin{propo}\label{Class-Stable}  
Let $K$ be an arbitrary field and $\cA$ a  $K$-algebra.  
Then $\cA$ is $\vartheta$-stable iff 
one of the following two statements holds:
\begin{enumerate}
  \item[$1)$] $\cA=K$;
  \item[$2)$] $K\simeq \bZ_2$ and 
  $\cA\simeq \bZ_2\dot{+}\bZ_2$.
  \end{enumerate}    
\end{propo}

\renewcommand{\theequation}{\thesection.\arabic{equation}}
\renewcommand{\therema}{\thesection.\arabic{rema}}
\setcounter{equation}{0}
\setcounter{rema}{0}

\section{\bf Mathieu Subspaces of Modules of Associative Algebras}
\label{S3} 

Let $R$ and $\cA$ be as before and $\cM$ an arbitrary 
(left) $\cA$-module. Recall that we have set   
in Section \ref{S1.2} the convention 
that all $\cA$-modules in this paper 
will be {\it left} $\cA$-modules unless stated otherwise. 
In particular, when we say {\it $\cA$ itself is an 
$\cA$-module}, we always mean that the left $\cA$-module 
 on the $R$-space $\cA$ with the action given 
by the algebra product of $\cA$ (from left).
%
%
%
%
%

Recall also that we have defined in 
Definition \ref{Def-GMS} the \vmss for  
the $\cA$-module $\cM$. 
But, for convenience we also fix the 
following terminologies.        

\begin{defi}\label{Def-Quasi-Stable}
With the notations and the setting as above, 
let $u\in \cM$ and $N$ a  $R$-subspace of 
$\cM$. We say that 

$1)$ $u$ is {\it $\vt$-stable}   
with respect to $N$ if $(N:u)$ 
is a $\vt$-ideal of $\cA$; 

$2)$ $u$ is {\it $\vt$-quasi-stable} 
with respect to $N$ if $(N:u)$ 
is a \vms of $\cA$. 
\end{defi}

By the definition above, the sets 
$\tau_\vt(N)$ and $\sigma_\vt(N)$ defined 
on page \pageref{TWO-SETS} in Section \ref{S1.2} 
are respectively the sets of {\it $\vt$-quasi-stable} 
elements and {\it $\vt$-stable} 
elements of the $R$-subspace $N$ 
of $\cM$.

The terminologies in Definition 
\ref{Def-Quasi-Stable} can be 
justified as follows for the case 
$\vt=${\it ``left"}. 

If $u$ is left stable with respect to $N$ and  
some $a\in \cA$ maps $u$ into the $R$-subspace 
$N$ of $\cM$, i.e., $au\in N$, then no element 
$b\in \cA$ can map $au$ outside of $N$. 
In other words, $au$ will stay inside 
$N$ ``forever" under the action of $\cA$.

If $u$ is left quasi-stable with respect to $N$ 
(equivalently, $N$ is a left Mathieu subspace 
of $\cM$ with respect to $u$) and for some 
$a\in \cA$, all the terms of the sequence 
$\{ a^m u\}$ lie inside $N$, 
then for any $b\in \cA$, all but finitely 
many terms of the sequence  
$\{b(a^mu)\}$ also lie inside $N$. 
In other words, no element $b\in \cA$ can map 
infinitely many terms of the sequence 
$\{a^mu\}$ outside the $R$-subspace $N$.  

The terminologies of the {\it right stableness}  
and {\it right quasi-stableness} 
can be interpreted parallelly  
but with slightly different 
meanings. 

From Definitions \ref{Def-GMS} and 
\ref{Def-Quasi-Stable}, it is easy to see that 
for every $R$-subspace $J$ of $\cA$, 
we have that $J$ is a \vms (resp., $\vt$-ideal) 
of the $R$-algebra $\cA$, iff $J$ as a  $R$-subspace of 
the $\cA$-module $\cA$ is $\vt$-quasi-stable (resp., $\vt$-stable) 
with respect to the identity element $1_\cA \in \cA$, 
iff $1_\cA\in \tau_\vt(N)$ (resp., $1_\cA\in \sigma_\vt(N)$). 

In this section, we mainly study some general properties of 
the sets $\tau_\vt(N)$ and $\sigma_\vt(N)$ 
of $\vt$-quasi-stable elements and $\vt$-stable 
elements, respectively, of $R$-subspaces 
$N$ of $\cA$-modules $\cM$. 
For a summary of the main results in this section, 
see the arrangement description for   
this section in Subsection \ref{S1.3}.

We start with the following two simple lemmas  
whose proofs are very straightforward 
and will be skipped here.  

\begin{lemma}\label{EasyLemma0}
For any $\cA$-module $\cM$, we have 
\begin{align}\label{EasyLemma0-e1}
\sigma_\vt(\cM) =\tau_\vt(\cM)=\cM.
\end{align}

Furthermore, when $\vt=${\it ``left"}, we have 
\begin{align}\label{EasyLemma0-e2}
\sigma_\vt(0) =\tau_\vt(0)=\cM.
\end{align}
\end{lemma}

Note that when $\vt\neq${\it ``left"}, 
Eq.\,(\ref{EasyLemma0-e2}) 
does not necessarily hold, e.g., 
see Proposition \ref{KnCase}, $ii)$ in later 
Section \ref{S5}. 

\begin{lemma}\label{EasyLemma1}
For each $R$-subspace $N$ 
of the $\cA$-module $\cM$, we have,    

$i)$ $\sigma_\vt(N)\subseteq \tau_\vt(N)$.

$ii)$ $0\in \sigma_\vt(N)$ and hence, 
$0\in  \tau_\vt(N)$.

$iii)$ $\sigma_\vt(N)$ and $\tau_\vt(N)$ are closed under the actions of the elements 
of $R^\times\cup\{0\}$. In particular, when the base ring 
$R$ is a field, $\sigma_\vt(N)$ and 
$\tau_\vt(N)$ are $R$-cones.
\end{lemma}

Actually, when $\vt=${\it ``left"}, 
$\sigma_\vt(N)$ is always closed 
under the action of $R$ even 
in the case that $R$ is not 
a field. More precisely, 
the following lemma also holds. 

\begin{lemma}\label{EasyLemma2}
Let $\vt=${\it ``left"} and $a\in \cA$. Then 
for each $R$-subspace $N\subseteq \cM$, 
we have $a\sigma_\vt(N)\subseteq \sigma_\vt(N)$.
\end{lemma}

\pf Let $u\in \sigma_\vt(N)$ (with $\vt=${\it ``left"}).
Then $(N:u)$ is a left ideal of $\cA$. 
Note that in general the following equation always holds:
\begin{align}\label{EasyLemma2-pe1}
(N:au)=\big( (N:u): a \big).
\end{align}
Furthermore, it is also well-known and easy to check that 
for every left ideal $J$ of $\cA$ and 
$b\in \cA$, $(J:b)$ is also a left 
ideal of $\cA$. Then it follows from this fact and 
Eq.\,(\ref{EasyLemma2-pe1}) that 
$(N:au)$ is also a left ideal of $\cA$. 
Therefore, we have $au\in \sigma_\vt(N)$, 
whence the lemma follows.  
\epfv

One remark on Lemma \ref{EasyLemma1}, $iii)$ 
and Lemma \ref{EasyLemma2} is that $\sigma_\vt(N)$ and $\tau_\vt(N)$ in general are not closed under the addition of $\cA$ 
(see the examples to be discussed in later 
Sections \ref{S4}--\ref{S6}).  

Next, we give some characterizations for the $\vt$-stable elements 
and $\vt$-quasi-stable elements of $R$-subspaces of $\cA$-modules, 
from which we also get some characterizations  
for \vmss of $\cA$-modules. 
 
\begin{lemma}
Let $u\in \cM$ and $N$ a  $R$-subspace 
of $\cM$. Set 
\begin{align}
a^{-1}N\!:=\{v\in \cM\,|\, av\in N\}. 
\label{a-Inv-N} 
\end{align}
Then the following statements holds: 
\begin{enumerate}
  \item[$i)$] $u$ is left stable with respect to $N$ iff for every $a\in \cA$, we have
\begin{align*}
(N:u)\subseteq (a^{-1}N: u),
\end{align*}

 \item[$ii)$]  $u$ is left quasi-stable with respect to $N$ iff for every $a\in \cA$, we have 
\begin{align*}
\sqrt{(N:u)}\subseteq  \sqrt{(a^{-1}N: u)}.
\end{align*}

 \item[$iii)$]  $u$ is right stable with respect to $N$ iff for every $a\in \cA$, we have
  \begin{align*}
(N:u)\subseteq  (N: au).
\end{align*}

 \item[$iv)$]  $u$ is right quasi-stable with respect to $N$ iff for every $a\in \cA$, we have 
\begin{align*}
\sqrt{(N:u)}\subseteq  \sqrt{(N: au)}.
\end{align*}

 \item[$v)$]  $u$ is $($two-sided$)$ stable with 
 respect to $N$ iff for every $a, b\in \cA$, we have
  \begin{align*}
(N:u)\subseteq  (a^{-1}N: bu).
\end{align*}

 \item[$vi)$]  $u$ is $($two-sided$)$ quasi-stable with respect to $N$ iff for every $a, b\in \cA$, we have 
\begin{align*}
\sqrt{(N:u)}\subseteq  \sqrt{(a^{-1}N: bu)}.
\end{align*}
\end{enumerate} 
\end{lemma}
\pf Note first that for any $a\in \cA$ and 
$u\in \cM$, by Eqs.\,(\ref{Def-Ann}), 
(\ref{a-Inv-J}) and (\ref{a-Inv-N}) 
we have the following equation: 
\begin{align}
(a^{-1}N:u)=a^{-1}(N:u). 
\end{align}
Then by the equation above and 
Eq.\,(\ref{EasyLemma2-pe1}), it is 
easy to see that the lemma follows 
immediately from Lemma \ref{CharsInRads} 
and Definition \ref{Def-Quasi-Stable}.
\epfv

The next lemma says that the operations of taking  
$\sigma_\vt(N)$ and $\tau_\vt (N)$ of $R$-subspaces $N$ 
of the $\cA$-module $\cM$ commute with the operation of 
taking the set $(S:u)$ of subsets $S\subseteq \cM$ 
with respect to elements $u\in \cM$.  

\begin{lemma}
Let $u\in\cM$ and $N$ be a  $R$-subspace 
of $\cM$. Then we have
\begin{align}
\big(\sigma_\vt(N):u\big) &=\sigma_\vt (N:u),
\label{uvt=vtu-e1} \\
\big(\tau_\vt(N) :u \big)   &=\tau_\vt(N:u).
\label{uvt=vtu-e2}
\end{align}
\end{lemma}
\pf For each $a\in \cA$,  
 we have that $a\in \big(\tau_\vt(N) :u \big)$, 
iff $au\in \tau_\vt(N)$, iff $(N:au)$ 
is a \vms of $\cA$, iff  
$\big((N:u):a\big)$ by 
Eq.\,(\ref{EasyLemma2-pe1}) is a \vms of $\cA$, 
iff $a\in \tau_\vt(N:u)$. 
Hence, Eq.\,(\ref{uvt=vtu-e2}) follows. 
Eq.\,(\ref{uvt=vtu-e1}) can be proved similarly. 
\epfv

The relations of the operations of taking  
$\sigma_\vt(N)$ and $\tau_\vt (N)$ 
on $R$-subspaces $N$ of $\cM$ with 
the operation of intersection is given 
by the following proposition.  

\begin{propo}\label{wIntersectn}
Let $N_i$ $(i\in I)$ be any collection 
of $R$-subspaces of $\cM$. 
Then  the following statements hold:  

$i)$ \begin{align}\label{wIntersectn-e1}
\bigcap_{i\in I} \, \sigma_\vt(N_i)  
\subseteq \sigma_\vt 
\Big(\bigcap_{i\in I} N_i \Big).
\end{align}

$ii)$ if the cardinal number $|I|<\infty$, then we have 
\begin{align}\label{wIntersectn-e2}
\bigcap_{i\in I} \, \tau_\vt(N_i)
\subseteq \tau_\vt \Big(\bigcap_{i\in I} N_i \Big).
\end{align} 
\end{propo}

\pf Note first that for each $u\in \cM$, 
by Eq.\,(\ref{Def-Ann}) we have
\begin{align}\label{wIntersectn-pe1}
\big( \big(\bigcap_{i\in I} N_i\big) : u\big)
=\bigcap_{i\in I} \, (N_i : u).
\end{align}
Then it is easy to see that Eq.\,(\ref{wIntersectn-e1}) follows from the equation above and the fact that the intersection of every collection of 
$\vt$-ideals of $\cA$ is also a $\vt$-ideal of $\cA$, and 
Eq.\,(\ref{wIntersectn-e2}) follows from the equation 
above and Lemma \ref{IntersectionLemma}. 
\epfv

\begin{corol}
Let $u\in \cM$ and $N_i$ 
$(1\le i\le k)$ \vmss of $\cM$ with 
respect to $u$. Then $\bigcap_{i=1}^k N_i$ is also a \vms of $\cM$ with respect to $u$.
\end{corol}

\begin{rmk}
Note that submodules and \vmss are closed 
under the operation of finite intersection 
$($see Lemma \ref{IntersectionLemma}$)$. 
By the corollary above, we see that this property 
is preserved under the generalizations $\mathcal G_i$ $(i=3, 4)$ in the diagram $(\ref{NotionCD})$ from the notions of submodules and \vmss, respectively, 
to the notion of \vmss of $\cA$-modules 
$\cM$. 
\end{rmk}

Next, for every $R$-subspace $N$ of $\cM$,  
we study the intersections of 
$\tau_\vt(N)$ and $\sigma_\vt(N)$ 
with $N$ itself. In contrast to 
the sets $\tau_\vt(N)$ and 
$\sigma_\vt(N)$ themselves, their 
intersections with $N$ are 
always $\cA$-submodules 
(see Theorem \ref{Max-submodule} below). 
Let's first prove the following lemma. 

\begin{lemma}\label{(u)-Lemma}
Let $u\in \cM$ and $N$ be a  $R$-subspace 
of $\cM$. Then we have

$i)$ $u\in N \cap \sigma_\vt(N)$ iff 
$\cA u \subseteq N$.

$ii)$ $u\in N \cap \tau_\vt(N)$ iff 
$\cA u \subseteq N$.
\end{lemma}

\pf Again, here we just give a proof for the case $\tau_\vt(N)$, 
i.e., for statement $ii)$. The proof of statement $i)$ 
is similar. 

$(\Rightarrow)$ Since $u\in \tau_\vt(N)$, 
we know that $(N:u)$ is a \vms 
of $\cA$. Since $u\in N$, we have 
$1_\cA \in (N:u)$. Then 
by Lemma \ref{ONE-Lemma}, 
$(N:u)=\cA$, whence $\cA u\subseteq N$.

$(\Leftarrow)$ Since $\cA u \subseteq N$, 
we have $u\in N$ and $(N:u)=\cA$. Since $\cA$  
obviously is a \vms of $\cA$, we have  
$u\in \tau_\vt(N)$, whence  
$u\in N \cap \tau_\vt(N)$.
\epfv

Recall that for each $R$-subspace $N$ of $\cM$, 
we have set in Subsection \ref{S1.3} $I_N$ to be  
the {\it maximum $\cA$-submodule of $N$}, 
i.e., the unique $\cA$-submodule of $\cM$  
which is maximum among all the 
$\cA$-submodules of $\cM$  
contained in $N$. It is easy to see that 
$I_N$ always exists and is actually 
the same as the sum 
of all the $\cA$-submodules of $\cM$  
contained in $N$. In particular, 
when $N$ itself is an $\cA$-submodule, 
we have $I_N=N$. 

\begin{theo}\label{Max-submodule}
For every $R$-subspace $N$ of $\cM$, we have 
\begin{align} \label{Max-submodule-e1}
I_N=N\cap\sigma_\vt(N)=N \cap\tau_\vt(N).
\end{align}
In particular, both $N\cap\tau_\vt(N)$ and 
$N \cap\sigma_\vt(N)$ are $\cA$-submodules  
of $\cM$.
\end{theo}

\pf Note that the second equality in 
Eq.\,(\ref{Max-submodule-e1}) follows directly 
from Lemma \ref{(u)-Lemma}. So it suffices  
to show the first equality in 
Eq.\,(\ref{Max-submodule-e1}). 

We first show that $N\cap \sigma_\vt(N)$ is closed 
under the addition of $\cM$. Let 
$u, v\in N\cap \sigma_\vt(N)$. Then  
by Lemma \ref{(u)-Lemma}, $i)$, we have 
$\cA u\subseteq N$ and $\cA v\subseteq N$. 
Hence, $\cA (u+v)\subseteq \cA u+\cA v\subseteq N$  
since $N$ is a  $R$-subspace of $\cM$.  
Then by Lemma \ref{(u)-Lemma}, $i)$ again, 
we have $u+v\in N\cap \sigma_\vt(N)$. 

Next, we show that $N\cap \sigma_\vt(N)$ 
is closed under the action of $\cA$. 
Let $a\in \cA$ and 
$u\in N\cap \sigma_\vt(N)$.  
By Lemma \ref{(u)-Lemma}, $i)$  
we have $\cA u \subseteq N$, whence  
$\cA (au)\subseteq \cA u \subseteq N$. 
Then by Lemma \ref{(u)-Lemma}, $i)$ again, 
we have $au \in N\cap \sigma_\vt(N)$.  

Therefore, $N\cap \sigma_\vt(N)$ is indeed an 
$\cA$-submodule of $\cM$ (contained in $N$). 
To show that $I_N=N\cap \sigma_\vt(N)$, 
it suffices to show that for each $\cA$-submodule 
$H\subseteq N$, we have $H\subseteq N \cap \sigma_\vt(N)$. 

Let $u\in H$. Since $H$ is an $\cA$-submodule 
of $\cA$, we have $\cA u\subseteq H\subseteq N$. 
Then by lemma \ref{(u)-Lemma}, $i)$,  
$u\in N\cap \sigma_\vt(N)$.  
Hence $H \subseteq N \cap \sigma_\vt(N)$.
\epfv

Next, we derive some consequences of 
Theorem \ref{Max-submodule}. 

\begin{corol}\label{Corol4Irred} 
For every specification of $\vt$ and an $\cA$-module 
$\cM$, the following statements are equivalent:

$1)$ $\cM$ is an irreducible 
$\cA$-module;

$2)$ for every proper $R$-subspace $N$ of $\cM$, 
we have $\sigma_\vt(N)\subseteq N^c\cup\{0\}$, 
where $N^c$ denotes the complement of $N$ in $\cM$; 

$3)$ for every proper $R$-subspace $N$ of $\cM$, 
we have $\tau_\vt(N)\subseteq N^c\cup\{0\}$.  
\end{corol}
\pf If statement $1)$ holds, i.e., $\cM$ is irreducible, 
then for every proper $R$-subspace $N$ of $\cM$, we have $I_N=0$. 
By Eq.\,(\ref{Max-submodule-e1}), we immediately 
have statements $2)$ and $3)$. 
Since by Lemma \ref{EasyLemma1}, $i)$  
$\sigma_\vt(N)\subseteq \tau_\vt(N)$ for every $R$-subspace 
$N$ of $\cM$, we also have $3)\Rightarrow 2)$. 
Therefore, it suffices to show $2)\Rightarrow 1)$.

Assume that $\cM$ is not irreducible. Let $N$ 
be a nonzero proper $\cA$-submodule of $\cM$. Then by 
Eq.\,(\ref{Max-submodule-e1}), we have 
$N \cap \sigma_\vt(N)=I_N=N\ne 0$, 
which contradicts statement $2)$.  
\epfv

%
%
%
%

\begin{corol}\label{SubmoduleCorol}
Let $N$ be a  $R$-subspace of $\cM$. 
Then the following statements are equivalent: 
\begin{enumerate}
\item[$1)$] $N\subseteq  \sigma_\vt(N)$; 

\item[$2)$] $N\subseteq  \tau_\vt(N)$;

\item[$3)$] $N$ is an $\cA$-submodule of $\cM$.
\end{enumerate}
\end{corol}
\pf Note that $N$ is an $\cA$-submodule of $\cM$ iff $I_N=N$. Then the corollary follows immediately from Eq.\,(\ref{Max-submodule-e1}).  
\epfv

Furthermore, for the case $\vt=${\it ``left"}, we can get $\sigma_\vt(N)$ and $\tau_\vt(N)$ explicitly under the equivalent conditions in the corollary above as follows. 

\begin{corol}\label{SubmoduleLemma}
Let $\vt=${\it ``left"} and $N$ be a  
$R$-subspace of $\cM$.
Then the following statements are equivalent: 
\begin{enumerate}
\item[$1)$] $N$ satisfies one of the equivalent 
statements $($with $\vt=${\it ``left"}$)$ in 
Corollary \ref{SubmoduleCorol};    

\item[$2)$] $\sigma_\vt(N)=\cM$;

\item[$3)$] $\tau_\vt(N)=\cM$. 
\end{enumerate} 
\end{corol}

\pf Since statements $2)$ and $3)$ 
obviously imply respectively statements $1)$ and $2)$ 
in Corollary \ref{SubmoduleCorol}, we have 
$2)\Rightarrow 1)$ and $3)\Rightarrow 1)$. 
Furthermore, by Lemma \ref{EasyLemma1}, $i)$ 
we also have $2)\Rightarrow 3)$. Therefore, 
it suffices to show $1)\Rightarrow 2)$. 

Assume that statement $3)$  
in Corollary \ref{SubmoduleCorol} holds, 
i.e., $N$ is a left $\cA$-submodule of $\cM$. 
Then it is easy to see that for any 
$u\in \cM$,  $(N:u)$ 
is also a left ideal of $\cA$, whence 
$u\in \sigma_\vt(N)$.  
Therefore, we have $\cM\subseteq \sigma_\vt(N)$, 
i.e., statement $2)$ holds.  
\epfv

\begin{rmk}\label{rmk-cG4} 
When $\vt=${\it ``left"}, from Corollaries \ref{SubmoduleCorol} 
and \ref{SubmoduleLemma} we see that 
for each $R$-subspace $N$ of the $\cA$-module 
$\cM$, $N$ is an $\cA$-submodule of $\cM$ iff $N$ is a left 
Mathieu subspace of $\cM$ with respect to all 
elements of $\cM$.  
Therefore, the notion of left Mathieu subspaces 
of $\cA$-modules in the sense above does 
generalize the notion of $\cA$-submodules.       
\end{rmk}

Next, we study some functorial properties of 
the sets of $\vt$-stable and $\vt$-quasi-stable 
elements of $R$-subspaces of $\cA$-modules.  

\begin{propo}\label{Inverse-vee-1}
Let $\cM$, $\cN$ be $\cA$-modules and 
$\phi: \cM\to \cN$ a homomorphism of 
$\cA$-modules. Then for every $R$-subspace $H$ 
of $\cN$, we have
\begin{align}
\phi^{-1}\big(\sigma_\vt(H)\big) 
&=\sigma_\vt\big(\phi^{-1}(H)\big), \label{Inverse-vee-1-e1} \\
\phi^{-1}\big(\tau_\vt(H)\big) 
&=\tau_\vt\big( \phi^{-1}(H) \big). \label{Inverse-vee-1-e2}
\end{align}
\end{propo}
\pf Note first that for each $u\in \cM$, it is easy to check 
that the following equation always holds:
\begin{align}
\big( H: \phi(u)\big)
=\big( \phi^{-1}(H): u \big). 
\label{Inverse-vee-1-pe1}
\end{align}

Then for each $u\in \cM$, 
we have that $u\in \tau_\vt\big( \phi^{-1}(H) \big)$,   
iff $\big( \phi^{-1}(H): u \big)$ 
is a \vms of $\cA$,   
iff $\big( H: \phi(u)\big)$ by 
Eq.\,(\ref{Inverse-vee-1-pe1}) is a \vms of $\cA$, 
iff $\phi(u)\in \tau_\vt(H)$, 
iff $u\in \phi^{-1}\big(\tau_\vt(H)\big)$. 
Hence, Eq.\,(\ref{Inverse-vee-1-e2}) follows. 
The proof of Eq.\,(\ref{Inverse-vee-1-e1}) 
is similar.
\epfv

Next we derive some consequences of 
Proposition \ref{Inverse-vee-1}. 

\begin{corol}
Let $a\in \cA$ be a central element of $\cA$, i.e., $a$ commutes with all elements of $\cA$. Then for every $R$-subspace 
$N$ of $\cM$, we have 
\begin{align}
a^{-1}\sigma_\vt(N)&=\sigma_\vt(a^{-1}N),\\
a^{-1}\tau_\vt(N) &=\tau_\vt(a^{-1}N). 
\end{align}
\end{corol}
\pf Let $\mu_a:\cM\to \cM$ be 
the $R$-linear map defined by the action 
of $a$ on $\cM$. Since $a$ is a central element 
of $\cA$, the map $\mu_a$ is an endomorphism 
of the $\cA$-module $\cM$. 
Then the corollary follows immediately from 
Proposition \ref{Inverse-vee-1} 
with $\phi=\mu_a$.
\epfv   

\begin{corol}\label{SubmoduleCase}
Let $V$ be an $\cA$-submodule of $\cM$ 
and $N$ a  $R$-subspace of $\cM$. 
Then we have
\begin{align}
V\cap \sigma_\vt(N)=\sigma_\vt(V\cap N), \label{SubmoduleCase-e1}\\
V\cap \tau_\vt(N)=\tau_\vt(V\cap N),\label{SubmoduleCase-e2}
\end{align}
where $\sigma_\vt(V\cap N)$ $($resp., $\tau_\vt(V\cap N)$$)$ 
is the set of $\vt$-stable $($resp., $\vt$-quasi-stable$)$ 
elements of $V\cap N$ as a $R$-subspace of the $\cA$-module $V$. 
\end{corol}

\pf Let $\iota: V\to \cM$ be the embedding of $V$ 
into $\cM$. Note that for every subset $S$ of $\cM$, 
we have $\iota^{-1}(S)=V\cap S$. Then 
the corollary follows immediately from 
this observation and 
Proposition \ref{Inverse-vee-1} 
with $\phi=\iota$. 
\epfv

\begin{corol}\label{GenerizePullBack}
In the same setting as 
in Proposition \ref{Inverse-vee-1}, let 
$u\in \cM$ and $v\in \cN$ such that 
$\phi(u)=v$. 
Then $H$ is \vms of $\cN$ 
with respect to $v\in \cN$ iff 
$\phi^{-1}(H)$ is \vms of $\cM$ 
with respect to $u\in \cM$. 
\end{corol}
\pf $H$ is a \vms of $\cN$ 
with respect to $v$, 
iff $\phi(u)=v \in \tau_\vt(H)$, iff  
$u\in  \phi^{-1}\big( \tau_\vt(H) \big)$, iff
$u\in \tau_\vt \big( \phi^{-1}(H) \big)$ (by 
Eq.\,(\ref{Inverse-vee-1-e2})),  
iff $\phi^{-1}(H)$ is a \vms of $\cM$ 
with respect to $u$. 
\epfv

\begin{rmk}\label{rmk-PullBack}
A fundamental property of submodules is that they are closed under the pull-backs of homomorphisms of 
modules. It is easy to see that 
Corollary \ref{GenerizePullBack} 
generalizes this property of submodules to \vmss 
of modules of associative algebras. Furthermore, 
Corollary \ref{GenerizePullBack} can also be viewed 
as a generalization of the similar property of \vmss 
of associative algebras given in 
Proposition \ref{Pull-Back} 
if we view the $R$-algebra $\cB$ in 
Proposition \ref{Pull-Back} as an $\cA$-module 
via the $R$-algebra homomorphism 
$\phi:\cA\to \cB$, and choose $u=1_\cA$ and $v=1_\cB$.    
\end{rmk}

Next, applying Proposition \ref{Inverse-vee-1} to  the quotient maps of associative algebras $\cA$, we get the following corollary, which can be viewed as a generalization of Proposition \ref{J-quotient}. 

\begin{corol}\label{Quo-Modules}
Let $\cM$ be an $\cA$-module and  
$V$ an $\cA$-submodule of $\cM$. Denote by $\pi: \cM\to \cM/V$ 
the quotient map. Then  
for every $u\in \cM$ and $R$-subspace  
$N$ of $\cM$ such that $V \subseteq N$, we have
\begin{align}
\sigma_\vt(N) &=\pi^{-1}\big(\sigma_\vt(N/V)\big), \label{Quo-Modules-e1} \\
\tau_\vt(N) &=\pi^{-1}\big(\tau_\vt(N/V)\big), \label{Quo-Modules-e2}\\
\pi\big(\sigma_\vt&(N)\big) =\sigma_\vt(N/V), \label{Quo-Modules-e1b} \\
\pi\big(\tau_\vt&(N)\big) =\tau_\vt(N/V). \label{Quo-Modules-e2b}
\end{align} 
\end{corol}
\pf Since $V\subseteq N$, we have 
$\pi^{-1}(N/V)=N$. Then by 
applying Proposition \ref{Inverse-vee-1} 
with $\phi=\pi$, 
we see that Eqs.\,(\ref{Quo-Modules-e1}) and (\ref{Quo-Modules-e2}) 
follow immediately from 
Eqs.\,(\ref{Inverse-vee-1-e1}) and (\ref{Inverse-vee-1-e2}), 
respectively. Furthermore, since $\pi$ is surjective, 
Eqs.\,(\ref{Quo-Modules-e1b}) and (\ref{Quo-Modules-e2b})
follow respectively from Eqs.\,(\ref{Quo-Modules-e1}) and (\ref{Quo-Modules-e2}).
\epfv

Finally, we conclude this section with the following variation of Proposition \ref{Inverse-vee-1}. 

\begin{propo}\label{Inverse-vee-2}
Let $\psi: \cA\to \cB$ be a  $R$-algebra  homomorphism. Then for every $R$-subspace 
$J$ of $\cB$, we have
\begin{align}
\psi^{-1}\big(\sigma_\vt(J)\big)
&\subseteq \sigma_\vt\big(\psi^{-1}(J)\big), \label{Inverse-vee-2-e1} \\
\psi^{-1}\big(\tau_\vt(J)\big)
&\subseteq \tau_\vt\big( \psi^{-1}(J) \big). \label{Inverse-vee-2-e2}
\end{align}
Furthermore, if $\psi$ is surjective, then the equalities in the both equations above hold. 
\end{propo}

\begin{rmk}
The reason that we do not always 
have equality in   
Eq.\,$($\ref{Inverse-vee-2-e2}$)$ 
$($in contrast to Eq.\,$($\ref{Inverse-vee-1-e2}$)$$)$ 
is because that the sets $\tau_\vt(J)$ and 
$\tau_\vt\big(\psi^{-1}(J)\big)$ of quasi-stable elements of $J$ 
and $\psi^{-1}(J)$, respectively, are   
defined with respect to the different 
algebras $\cB$ and $\cA$. 
This is also the case for 
Eq.\,$($\ref{Inverse-vee-2-e1}$)$. 
\end{rmk} 

\underline{\it Proof of 
Proposition \ref{Inverse-vee-2}}\,: 
Note first that for each $a\in\cA$, it is easy  
to check directly that the following equation holds: 
\begin{align}
\psi^{-1}\big(J:\psi(a)\big)=
\big(\psi^{-1}(J): a\big).
\label{Inverse-vee-2-pe1} 
\end{align}

Note also that the sets $\big(J:\psi(a)\big)$ and 
$\big(\psi^{-1}(J): a\big)$ in the equation above 
are defined with respect 
to the (different) $R$-algebras $\cB$ and 
$\cA$, respectively.

Now, for each  
$a\in \psi^{-1}\big(\tau_\vt(J)\big)$,  
we have $\psi(a)\in \tau_\vt(J)$, 
whence $(J:\psi(a))$ 
is a \vms of the $R$-algebra $\cB$.
Then by Eq.\,(\ref{Inverse-vee-2-pe1}) and 
Proposition \ref{Pull-Back}, 
$\big(\psi^{-1}(J): a\big)$ is a \vms of $\cA$. 
Hence, we have $a\in \tau_\vt\big( \psi^{-1}(J)\big)$, 
whence Eq.\,(\ref{Inverse-vee-2-e2}) follows.  
Eq.\,(\ref{Inverse-vee-2-e1}) can be proved 
similarly by using Eq.\,(\ref{Inverse-vee-2-pe1}) 
and the fact that $\vt$-ideals are closed 
under the pulling-back by $R$-algebra 
homomorphisms.

Next we assume that $\psi$ is surjective and show 
that the equality in Eq.\,(\ref{Inverse-vee-2-e2}) 
does hold. The equality in 
Eq.\,(\ref{Inverse-vee-2-e1}) 
can be proved similarly. 

Let $a \in \tau_\vt \big( \psi^{-1}(J) \big)$. 
Then $\big(\psi^{-1}(J): a\big)$ is a \vms of $\cA$. 
We need to show $a \in \psi^{-1}\big (\tau_\vt(J) \big)$  
or equivalently, $(J:\psi(a))$ is a \vms of $\cB$.

Let $I$ be the kernel of the $R$-homomorphism 
$\psi$. Then $I$ is an ideal of $\cA$. By   
identifying the $R$-algebra $\cB$ with 
the quotient $R$-algebra $\cA/I$ via 
the induced $R$-algebra isomorphism 
$\bar \psi: \cA/I\overset{\sim}{\to} \cB$,  
we may assume that $\cB=\cA/I$ and 
the $R$-algebra homomorphism 
$\psi$ is the quotient map 
$\pi: \cA \to \cA/I$.    

With the setting above, it is 
easy to check that we have 
\begin{align}\label{Inverse-vee-2-pe2} 
I\subseteq \big(\psi^{-1}(J):a\big). 
\end{align}

Combining the equation above with 
Eq.\,(\ref{Inverse-vee-2-pe1}), we also have 
\begin{align}\label{Inverse-vee-2-pe3} 
\big(\psi^{-1}(J):a\big)/I=\big(J: \psi(a)\big ). 
\end{align}

Then with Eqs.\,(\ref{Inverse-vee-2-pe2}) 
and (\ref{Inverse-vee-2-pe3}) above, applying  Proposition \ref{J-quotient} to the \vms 
$\big(\psi^{-1}(J):a\big)$ of $\cA$, we see that 
$(J:\psi(a))$ is indeed a \vms of $\cB$. 
\epfv

%
%

\renewcommand{\theequation}{\thesection.\arabic{equation}}
\renewcommand{\therema}{\thesection.\arabic{rema}}
\setcounter{equation}{0}
\setcounter{rema}{0}

\section{\bf Quasi-Stable Modules and Quasi-Stable Algebras}\label{S4}

Let $R$, $\cA$ and $\cM$ be fixed as before. 
We first introduce the following notions. 

\begin{defi}\label{Def-QSM}
An $\cA$-module $\cM$ is said to be 
{\it $\vartheta$-quasi-stable} 
$($resp., $\vartheta$-stable$)$ if 
for every $R$-subspace $N$ of $\cM$, 
we have $N^c\subseteq  \tau_\vt(N)$ 
$($resp., $N^c\subseteq  \sigma_\vt(N)$$)$, 
where $N^c$ denotes the complement 
of the subset $N$ in $\cM$. 
\end{defi}



Note that by Lemma \ref{EasyLemma1}, $i)$, we see that   
 the $\vt$-stableness implies 
$\vt$-quasi-stableness.
 
From Definition \ref{Def-QSM}, Theorem \ref{Max-submodule} 
and also the  observation above, it is easy to see that 
the following lemma holds. 

\begin{lemma}\label{DirectLemma} 
$i)$\, If $\cM$ is a $\vt$-quasi-stable $\cA$-module, 
then for any  $R$-subspace 
$N$ of $\cM$, we have
\begin{align}\label{DirectLemma-e1} 
\tau_\vt (N)=I_N \cup N^c.
\end{align}  

$ii)$\, If $\cM$ is a $\vt$-stable $\cA$-module, 
then for any 
$R$-subspace $N$ of $\cM$, we have
\begin{align}\label{DirectLemma-e2} 
\sigma_\vt (N)=\tau_\vt(N)=I_N \cup N^c.
\end{align}  
\end{lemma} 

When the $R$-algebra $\cA$ itself 
is viewed as a left $\cA$-module (in the canonical way), 
then the next lemma shows that Definitions \ref{q-StaAlg} 
and \ref{Def-QSM} actually coincide.

\begin{lemma}\label{qStableEquiv-Lemma}
For every $R$-algebra $\cA$, it is {\it $\vt$-quasi-stable} $($resp., $\vt$-stable$)$ as a  $R$-algebra iff it is {\it $\vt$-quasi-stable} 
$($resp., $\vt$-stable$)$ as a left $\cA$-module. 
\end{lemma}

\pf We just give a proof here for the $\vt$-quasi-stable case. The $\vt$-stable case can be proved similarly. 

$(\Rightarrow)$ Let $J$ be a  $R$-subspace 
of $\cA$ and $a\in J^c$. Then  
$1\not \in (J:a)$ and the $R$-subspace 
$(J:a)$ is a \vms of $\cA$, for 
$\cA$ is a $\vt$-quasi-stable  
$R$-algebra, whence $a\in \tau_\vt(J)$.
Therefore, we have $J^c\subseteq \tau_\vt(J)$.
Hence $\cA$ as a left 
$\cA$-module is also    
$\vt$-quasi-stable.

$(\Leftarrow)$ Let $J$ be a  $R$-subspace 
of $\cA$ with $1 \not \in J$. 
Since $\cA$ as a left $\cA$-module 
is $\vt$-quasi-stable, we have 
$J^c\subseteq  \tau_\vt(J)$. 
Since $1\in J^c$, we have $1\in \tau_\vt(J)$, 
whence $J=(J:1)$ is a \vms of $\cA$. Therefore, 
$\cA$ as a  $R$-algebra is also 
$\vt$-quasi-stable.
\epfv

Actually, the following more general result 
also holds.

\begin{propo}\label{qStableEquiv}
For every $R$-algebra $\cA$, the following two 
statements are equivalent to each other.

$1)$ $\cA$ as a  $R$-algebra is {\it $\vt$-quasi-stable} 
$($resp., $\vt$-stable$)$. 

$2)$ Every left $\cA$-module $\cM$ is {\it $\vt$-quasi-stable} 
$($resp., $\vt$-stable$)$.
\end{propo}
\pf $2)\Rightarrow 1)$ follows from Lemma 
\ref{qStableEquiv-Lemma}. 
$1)\Rightarrow 2)$ can be shown by 
a similar argument as in the proof of the 
$(\Rightarrow)$ part of Lemma 
\ref{qStableEquiv-Lemma}, 
at least for the $\vt$-quasi-stable 
case. The $\vt$-stable case 
can also be proved similarly. 
\epfv

Combining the proposition above with 
Lemma \ref{DirectLemma},   
Proposition \ref{q-Quasi/R},  
Theorem \ref{Class-qStable} and 
Proposition \ref{AlgLocal}, 
we immediately 
have the following corollary.  

\begin{corol}\label{Tau4qStable} 
Let $\cA$ be a  $R$-algebra as in 
Proposition \ref{q-Quasi/R} 
or a $K$-algebra satisfying one of the  
equivalent statements 
in Proposition \ref{AlgLocal}, 
or $\cA\simeq K\dot{+}K$. 
Then all left $\cA$-modules 
$\cM$ are $\vt$-quasi-stable. Consequently, 
for any $R$-subspace $N$ of $\cM$, we have 
$\tau_\vt(N)=I_N\cup N^c$, which is actually 
independent on the specifications of $\vt$.    
\end{corol}


\begin{corol}\label{DivisionCase}
Let $\cA$ be an algebraic field extension of a field $K$ or more generally, an algebraic division algebra over $K$. Then for every 
$K$-subspace $J$ of $\cA$, the following two statements hold. 

$i)$ If $J=0$ or $\cA$, then we have  
\begin{align}\label{DivisionCase-e1}
\sigma_\vt(J)=\tau_\vt (J)=\cA.
\end{align}  

$ii)$ If $J$ is nonzero and proper, then we have 
\begin{align}
\sigma_\vt(J)&=\{0\}, \label{DivisionCase-e2}\\
\tau_\vt (J) &= J^c\cup\{0\}. 
\label{DivisionCase-e3}
\end{align} 
\end{corol} 
\pf Note first that for every $a\in \cA$, 
we have 
\begin{align}\label{DivisionCase-pe1}
(\cA:a)=\cA. 
\end{align} 
Since $\cA$ has no zero-divisors, 
we also have 
\begin{align}\label{DivisionCase-pe2}
(0:a)=
\begin{cases}
 \cA &\mbox{ if } a=0;\\
0 &\mbox{ if } a\ne 0. 
\end{cases}
\end{align} 
Then it is easy to see that 
Eq.\,(\ref{DivisionCase-e1}) follows immediately 
from the two equations above.

To show Eq.\,(\ref{DivisionCase-e2}), 
assume that $J$ is nonzero and proper. 
Note first that for each  
$0\ne a\in \cA$ and $0\ne b\in J$, 
we have $0\ne ba^{-1} \in (J:a)$, whence  
$(J:a)\ne 0$. Moreover, 
$(J:a)\ne \cA$ either, 
for otherwise we would have 
$\cA=\cA a\subseteq J$, 
which contradicts our assumption 
that $J$ is proper.  Therefore, $(J:a)$ 
is also a nonzero proper $K$-subspace 
of $\cA$. 

On the other hand, since $\cA$ has no nonzero proper 
$\vt$-ideals (for every nonzero element of $\cA$ is invertible),  
the $K$-subspace $(J:a)$ cannot be a $\vt$-ideal of $\cA$. 
Hence we have $a \not \in \sigma_\vt(J)$ 
for all $0\ne a\in \cA$, whence 
Eq.\,(\ref{DivisionCase-e2}) 
follows.    

To show Eq.\,(\ref{DivisionCase-e3}), note first that 
by Theorem \ref{Class-qStable} and Proposition \ref{AlgLocal},  
$\cA$ is a $\vt$-quasi-stable $K$-algebra. 
Then by Lemma \ref{qStableEquiv-Lemma} or 
Proposition \ref{qStableEquiv}, 
$\cA$ as a left $\cA$-module is also 
$\vt$-quasi-stable. Therefore,  
by Eq.\,(\ref{DirectLemma-e1}) we have 
$\tau_\vt (J)=I_J \cup J^c$, where 
$I_J$ is the (unique) left ideal of $\cA$ 
which is maximum among all the left ideals 
of $\cA$ contained in $J$.

On the other hand, since $\cA$ has no nonzero  
proper left ideals (as already pointed out above), 
we have $I_J=0$, whence Eq.\,(\ref{DivisionCase-e3}) 
follows.  
\epfv

\begin{rmk}\label{IrredWn0tau} 
Let $\cA$ be as in 
Corollary \ref{DivisionCase}. Then 
$\cA$ as a left $\cA$-module is 
irreducible, since every nonzero element 
of $\cA$ is invertible. Therefore,  
Corollary \ref{DivisionCase} provides 
a family of irreducible $\cA$-modules 
$\cM$ such that $\tau_\vt(J)\ne 0$ 
for all proper $K$-subspaces $J$ of 
$\cM$. 
\end{rmk}

Finally, we conclude this section with the following 
open problem which, we believe, is worthy to be 
much further investigated.  
\begin{prob}
Classify all $\vt$-stable or $\vt$-quasi-stable modules for some 
``nice"  associative algebras, say, for semi-simple algebras 
or Noetherian algebras. 
\end{prob}

\renewcommand{\theequation}{\thesection.\arabic{equation}}
\renewcommand{\therema}{\thesection.\arabic{rema}}
\setcounter{equation}{0}
\setcounter{rema}{0}

\section{\bf Two Cases for Modules of Matrix Algebras 
over Fields}\label{S5}

In this section, we fix an arbitrary field $K$ 
and let $M_n(K)$ $(n\ge 1)$ denote  
the matrix algebra of $n\times n$ matrices 
with entries in $K$. 

We first consider the sets $\sigma_\vt(N)$ 
and $\tau_\vt(N)$ for $K$-subspaces $N$ of 
the $M_n(K)$-module $K^n$ 
(with the standard left action).   

\begin{propo}\label{KnCase}
Let $n\ge 2$ and $N$ a  $K$-subspace 
of the $M_n(K)$-module $K^n$. Then 
the following statements hold.

$i)$ If $N=K^n$, we have 
\begin{align}\label{KnCase-e1}
\sigma_\vt(N)=\tau_\vt(N)=K^n.
\end{align}

$ii)$ If $N=0$, we have 
\begin{align}\label{KnCase-e2}
\sigma_\vt(N)=\tau_\vt(N)=
\begin{cases}
K^n &\mbox{ if } \vt=``left";\\
0 &\mbox{ otherwise.} 
\end{cases} 
\end{align}

$iii)$ If $N$ is nonzero and proper, we have 
\begin{align}\label{KnCase-e3} 
\sigma_\vt(N)=\tau_\vt(N)=0.
\end{align}
\end{propo}

Note that the case $n=1$ has been covered by Corollary \ref{DivisionCase} in the previous section. More precisely, in this case Eq.\,(\ref{KnCase-e1}) still holds but not 
Eq.\,(\ref{KnCase-e2}). Instead, 
we have $\sigma_\vt(0)=\tau_\vt(0)=K$ for every specification of $\vt$. 

In order to prove the proposition above, we first need to show the following lemma.

\begin{lemma}\label{KnCase-Lemma} 
Let $n\ge 2$ and $u, v\in K^n$ such that 
$u$ and $v$ are $K$-linearly independent 
with each other. 
Then the following two statements hold:

$i)$ there exists an idempotent 
$E_1 \in M_n(K)$ 
such that $E_1 u=v$ and $E_1 v=v$; 


$ii)$ there exists an  
idempotent $E_2\in M_n(K)$ 
such that $E_2 u=0$ and $E_2 v=u+v$.
\end{lemma}

\pf We identify $M_n(K)$ with the $K$-algebra 
${\rm End}_K(K^n)$ of $K$-linear endomorphisms  
of $K^n$ via the standard basis of $K^n$.
 
Let $\{v_1,v_2, ..., v_n\}$ be a  $K$-linear basis of $K^n$ such that $u=v_1$ and $v=v_2$. Let $E_1$ be the $K$-linear map such that 
$E_1 v_1=v_2$ and $E_1 v_i=v_i$ for all 
$2\le i\le n$. Then it is easy to check that 
$E_1^2=E_1$, whence $i)$ follows. 

To show $ii)$, let $E_2$ be the $K$-linear map 
such that $E_2 v_1=0$;  
$E_2 v_2=v_1+v_2$ and $E_2 v_i=v_i$ 
for all $3\le i\le n$. 
Then it is easy to check that $E_2^2=E_2$,  
whence $ii)$ follows.   
\epfv

\underline{\it Proof of 
Proposition \ref{KnCase}}\,: 
Note first that Eq.\,(\ref{KnCase-e1}) and also 
the case $\vt=${\it ``left"} of Eq.\,(\ref{KnCase-e2}) 
follow directly from Lemma \ref{EasyLemma0}. 

To show the case $\vt=${\it ``right"} of  
Eq.\,(\ref{KnCase-e2}), by Lemma \ref{EasyLemma1}, $i)$, 
it suffices to show $\tau_\vt(0)=0$ for $\vt=${\it ``right"}. 
More explicitly, it suffices to show that 
the annihilator $(0:u)$ for each 
$0\ne u\in K^n$ is not a 
right Mathieu subspace of $M_n(K)$. 

We fix an arbitrary $0\ne u\in K^n$ and pick up a nonzero 
$v\in K^n$ such that $u$ and $v$ are 
$K$-linearly independent. Applying Lemma \ref{KnCase-Lemma} to the vectors $u$ and $v$, and letting $E_i$ $(i=1, 2)$ be the idempotents in the same lemma, we have $E_2\in (0:u)$ 
but $E_2E_1\not \in (0:u)$, for $E_2E_1u=E_2v=u+v\ne 0$. Hence, the right ideal of $M_n(K)$ generated by the idempotent $E_2$ is not contained in $(0:u)$. Then by 
Theorem \ref{CharByIdem}, $(0:u)$ is not a right Mathieu subspace of $M_n(K)$. Therefore, the case $\vt=${\it ``right"} of Eq.\,(\ref{KnCase-e2}) holds, whence so do the cases  $\vt=${\it``pre-two-sided"} and 
$\vt=${\it``two-sided"}.  

To show Eq.\,(\ref{KnCase-e3}), by Lemma \ref{EasyLemma1}, $i)$, it suffices to show $\tau_\vt(N)=0$. Moreover, 
since $K^n$ is an irreducible $M_n(K)$-module, by Corollary \ref{Corol4Irred}, we have 
$\tau_\vt(N) \subseteq N^c\cup\{0\}$. Therefore, 
it suffices to show 
\begin{align}\label{KnCase-pe1}
N^c\cap \tau_\vt(N)=\emptyset. 
\end{align}

Now, let $u\in N^c$ and choose any $0\ne v\in N$. 
Note that $u$ and $v$ are $K$-linearly 
independent. Applying Lemma \ref{KnCase-Lemma} to 
$u$ and $v$, and letting $E_i$ $(i=1, 2)$ 
be the idempotents in the same lemma, 
we have $E_1, E_2\in (N:u)$.  

Let $A\in M_n(K)$ such that $Au=v$ and $Av=u$   
(Note that by identifying $M_n(K)$ with ${\rm End}_K(K^n)$,  
as in the proof of Lemma \ref{KnCase-Lemma}, it is easy to see 
that such a matrix $A$ does exist). 
Then we have  
$AE_1 u=u\not \in N$. Hence, 
$AE_1\not \in (N:u)$ 
and by Theorem \ref{CharByIdem},  
$(N:u)$ is not a left 
Mathieu subspace of $M_n(K)$.
On the other hand, we also have 
$E_2 A u=E_2v =u+v\not \in N$. 
Hence, $E_2A \not \in (N:u)$, 
and by Theorem \ref{CharByIdem},  
$(N:u)$ is not a right 
Mathieu subspace of $M_n(K)$ 
either. Consequently, 
$u\not \in \tau_\vt(N)$ when 
$\vt=${\it ``left"} or {\it ``right"}, 
and hence the same when $\vt=${\it``pre-two-sided"} and 
$\vt=${\it``two-sided"}. Therefore, Eq.\,(\ref{KnCase-pe1}) 
does hold for all specifications of $\vt$. 
\epfv

Next, we consider the sets $\sigma_\vt(J)$ 
and $\tau_\vt(J)$ for co-dimension one $K$-subspaces $J$ 
of $M_n(K)$ as a left $M_n(K)$-module. 
First, let's fix the following notations.

We denote by $I_n$ the identity $n\times n$ matrix  
and ${\rm Tr\,}: M_n(K)\to K$ the {\it trace} function of 
$M_n(K)$. Furthermore, we set 
\begin{align}\label{Def-VX}
H_X\!:=\{A\in M_n(K)\,|\, {\rm Tr\,} (AX) =0\}.
\end{align} 
For simplicity, we also denote $H_{I_n}$ by 
$H$, i.e., $H$ is the co-dimension one $K$-subspace of 
the trace-zero matrices of $M_n(K)$. 
  
Note that the trace function ${\rm Tr\,}:M_n(K) \to K$ induces the following non-singular $K$-bilinear form of 
$M_n(K)$: 
\begin{align}
(\cdot, \cdot): M_n(K)\times M_n(K) &\to \quad  K  
\label{T-Pairing} \\
(Y,\quad X)\quad \quad &\to \,  {\rm Tr\,} (YX). \nno 
\end{align}

Let $X, Y\in M_n(K)$. We denote by $X\sim Y$ 
if $X=\alpha Y$ for some $\alpha\in K^\times$. 
Note that by Eq.\,(\ref{Def-VX}) and the 
non-singularity of the $K$-bilinear form in 
Eq.\,(\ref{T-Pairing}), it is easy to check that 
\begin{align}
H_X=H_Y \, \Leftrightarrow \, X\sim Y. \label{NX=NY} 
\end{align}
In particular, we have  
\begin{align}
H_X=M_n(K)  \, &\Leftrightarrow \, X=0, \label{HX=MK} \\
H_X=H \, &\Leftrightarrow \, X \sim I_n. \label{NX=H} 
\end{align}

Furthermore, it is also easy to check 
(e.g., see Lemma $5.2$ in \cite{MS}) that 
every co-dimension one $K$-subspace of $M_n(K)$ 
has the form $H_X$ for some $0\ne X\in M_n(K)$   
and the matrix $X$ by Eq.\,(\ref{NX=NY}) is 
unique up to nonzero scalar multiplications. 

With the notations fixed above, we can state our second main result of this section as follows.

\begin{propo}\label{V-CoD1}
Let $K$ be a field and $n\ge 2$. Then for each  
$X\in M_n(K)$, the following statements hold. 

$i)$ If $char.\,K=0$ or $char.\,K=p>n$, we have 
\begin{align}
\sigma_\vt(H_X)&=(0:X)=\{Y\in M_n(K)\,|\, YX=0\}. 
\label{V-CoD1-e1}\\
\tau_\vt(H_X)
&=\{Y\in M_n(K)\,|\, YX=0\, \mbox{ or }\, YX\sim I_n\}.
\label{V-CoD1-e2}
\end{align}

$ii)$ If $char.\,K=p>0$ and $p\le n$, then 
\begin{align}
\sigma_\vt(H_X)=\tau_\vt(H_X)=(0:X).\label{V-CoD1-e3}
\end{align}
\end{propo}

Note that the proposition above requires $n\ge 2$. The case 
$n=1$ has been covered by Corollary \ref{DivisionCase} 
in Section \ref{S4}. Note also that for all $n\ge 1$ and 
$X\in M_n(K)$, 
by Corollary \ref{DivisionCase} and the proposition above 
the sets $\sigma_\vt(H_X)$ and $\tau_\vt(H_X)$ are 
actually independent on the specifications 
of $\vt$. 

In order to prove the proposition above, we need the following theorem proved in \cite{MS} on the co-dimension one \vmss of $M_n(K)$.

\begin{theo}\label{Co-D1}
Let $K$ be a field and $n\ge 1$. Then for every fixed 
specification of $\vt$, the following two statements hold.

$i)$ If $char.\,K=0$ or $char.\,K=p>n$, then 
$H$ is the only co-dimension one 
\vms of $M_n(K)$.

$ii)$ If $char.\,K=p>0$ and $p\le n$, then 
$M_n(K)$ has no co-dimension one 
$\vt$-Mathieu subspace.
\end{theo}

\underline{\it Proof of Proposition \ref{V-CoD1}}\,: 
Let $Y\in M_n(K)$. Then by the definition 
of $H_X$ in Eq.\,(\ref{Def-VX}) it is easy 
to see that the following equation holds: 
\begin{align}\label{XY=YX}
(H_X:Y)=H_{YX}. 
\end{align}

Assume first $YX=0$, i.e., $Y\in (0:X)$.   
Then by the equation above and 
Eq.\,(\ref{HX=MK}), 
we have $(H_X:Y)=M_n(K)$. Since 
$M_n(K)$ is obviously a $\vt$-ideal of 
$M_n(K)$, we have $Y\in \sigma_\vt(H_X)$.  
Consequently, we have $(0:X) \subseteq \sigma_\vt(H_X)$.    
By Lemma \ref{EasyLemma1}, $i)$, we also have 
\begin{align}\label{V-CoD1-pe2}
(0:X) \subseteq \sigma_\vt(H_X) \subseteq \tau_\vt(H_X). 
\end{align} 

Assume $YX\ne 0$ and $YX \not\sim I_n$. 
Then by Eq.\,(\ref{NX=H}), we have 
$H_{YX}\ne H$, and by  
Theorem \ref{Co-D1} and 
Eq.\,(\ref{XY=YX}), $(H_X:Y)$ is not a 
\vms of $M_n(K)$ and hence, 
not a $\vt$-ideal of $M_n(K)$ either. 
Therefore, in this case we have 
\begin{align}\label{V-CoD1-pe3}
Y\not \in \tau_\vt(H_X) \, 
\mbox{ and } \, 
Y\not \in \sigma_\vt(H_X). 
\end{align}

Now, assume $YX \sim I_n$. 
Then by Eqs.\,(\ref{NX=H}) and 
(\ref{XY=YX}), we have $(H_X:Y)=H$. 
Furthermore, by Theorem \ref{Co-D1}  
we see that $(H_X:Y)$ is a \vms of 
$M_n(K)$, i.e., $Y\in \tau_\vt(H_X)$, 
iff the condition of the statement 
$i)$ holds. 

From the arguments above, it is easy to see that 
Eq.\,(\ref{V-CoD1-e3}) for $\tau_\vt(H_X)$ and also 
Eq.\,(\ref{V-CoD1-e2}) indeed hold, i.e., 
the proposition holds for the set 
$\tau_\vt(H_X)$. Furthermore, by Eq.\,(\ref{V-CoD1-pe2}) 
and Eq.\,(\ref{V-CoD1-e3}) for $\tau_\vt(H_X)$, 
we also have  Eq.\,(\ref{V-CoD1-e3}) for $\sigma_\vt(H_X)$.  

Finally, to show Eq.\,(\ref{V-CoD1-e1}), 
by Eqs.\,(\ref{XY=YX})-(\ref{V-CoD1-pe3}) 
it suffices to show that the co-dimension 
one $K$-subspace $H$ is not a left or right ideal 
of $M_n(K)$. But this can be easily checked under 
the condition $n\ge 2$. Therefore, 
the proposition follows. 
\epfv 

One consequence of 
Proposition \ref{V-CoD1} is the 
following corollary, which can also be proved 
directly by using the non-singularity 
of the paring in Eq.\,(\ref{T-Pairing}). 

\begin{corol}\label{V-CoD1-Corol}
Let $K$ be a field and $n\ge 1$. For each 
$X\in M_n(K)$, denote by $I_{X}$ $($resp., $J_X$$)$ 
the unique left $($resp., right$)$ ideal of 
$M_n(K)$ which is maximum among all the left 
$($resp., right$)$ ideals of $M_n(K)$ 
contained in $H_X$.  Then we have 
\begin{align}
I_X&=(0:X)=\{Y\in M_n(K)\,|\, YX=0\},
\label{V-CoD1-Corol-e1}\\
J_X&=\{Y\in M_n(K)\,|\, XY=0\}. \label{V-CoD1-Corol-e2}
\end{align}
%
%
\end{corol}

\pf The case $n=1$ can be checked easily. So we assume 
$n\ge 2$.

Note that by Proposition \ref{V-CoD1} 
we have 
\begin{align}\label{V-CoD1-Corol-pe1}
\sigma_\vt(H_X)=\{Y\in M_n(K)\,|\, YX=0\}.
\end{align}
Then viewing $M_n(K)$ as a left 
$M_n(K)$-module and applying 
Theorem \ref{Max-submodule}, we see 
that Eq.\,(\ref{V-CoD1-Corol-e1}) follows 
from Eq.\,(\ref{V-CoD1-Corol-pe1}) 
for the case $\vt=${\it ``left"}.

To show Eq.\,(\ref{V-CoD1-Corol-e2}), we first consider 
the {\it transpose} map 
\begin{align*}
t: M_n(K)&\to M_n(K) \\
Y\, &\to \,Y^t, 
\end{align*}
where $Y^t$ denotes the  {\it transpose} 
of the matrix $Y$. 

Note that the map $t$ is an {\it anti-involution} of  
$K$-algebra $M_n(K)$, i.e., $t^2={\rm id}$,  
the identity map of $M_n(K)$, and $t(XY)=t(Y)t(X)$ 
for all $X, Y\in M_n(K)$. Therefore, $t$ induces an one-to-one 
correspondence between the set of 
all left ideals $J$ of $M_n(K)$ and the set 
of all right ideals of $M_n(K)$ 
via $J \leftrightarrow t(J)$.

Furthermore, for all $X, Y\in M_n(K)$, we also have 
\begin{align}\label{V-CoD1-Corol-pe2}
{\rm Tr\,}(YX)={\rm Tr\,}\big((YX)^t\big)
={\rm Tr\,}(X^tY^t)={\rm Tr\,}(Y^tX^t).
\end{align}

From Eq.\,(\ref{Def-VX}) and the equation above, 
it is easy to see that for each $X\in M_n(K)$, 
we have 
\begin{align}\label{V-CoD1-Corol-pe3}
t(H_X)=H_{X^t}.
\end{align}

Now, by Eq.\,(\ref{V-CoD1-Corol-e1}) with 
$X$ replaced by $X^t$, we have  
\begin{align}
I_{X^t}&=(0:X^t)=\{Y\in M_n(K)\,|\, YX^t=0\}.
\end{align}

Then by the last two equations above, we have 
\allowdisplaybreaks{
\begin{align*}
J_X&=t^{-1}(I_{X^t})=t\,(I_{X^t})
=\{ Y^t\in M_n(K) \,|\, YX^t=0\} \\
&=\{ Y \in M_n(K)\,|\, Y^tX^t=0\}\\
&=\{ Y \in M_n(K)\,|\, (XY)^t=0\} \\
&=\{ Y \in M_n(K)\,|\, XY=0\}.
\end{align*} }
Hence, we get Eq.\,(\ref{V-CoD1-Corol-e2}). 
\epfv


\renewcommand{\theequation}{\thesection.\arabic{equation}}
\renewcommand{\therema}{\thesection.\arabic{rema}}
\setcounter{equation}{0}
\setcounter{rema}{0}

\section{\bf Two Cases for Modules of Polynomial Algebras over Fields}
\label{S6}
Let $K$ be a field and $z=(z_1, z_2, ..., z_n)$ $n$ commutative free variables. In this section, 
we consider the sets $\sigma_\vt(N)$ and $\tau_\vt(N)$ of the stable elements and quasi-stable elements, respectively, of two families of co-dimension one $K$-subspaces $N$ 
of the polynomial algebra $K[z]$ (as a  $K[z]$-module).

Note that $K[z]$ is commutative. The sets 
$\sigma_\vt(N)$ and $\tau_\vt(N)$ are actually 
independent on the specifications of $\vt$. 
Therefore, we simply denote them by 
$\sigma(N)$ and $\tau(N)$, respectively. 
Furthermore, we also need to 
fix the following notations. 

For any $\ell\ge 1$ and 
$\alpha=(\alpha_1, \alpha_2, ..., 
\alpha_\ell) \in K^\ell$, we denote by 
$S_\alpha$ the set of $1\le i\le \ell$ 
such that $\alpha_i\ne 0$. We let    
$\Omega_\ell$ be the set 
of $\alpha=(\alpha_1, \alpha_2, ...,  
\alpha_\ell)\in K^\ell$, which satisfies  
the following property:  

\vspace{2mm}
$(\ast)$ \quad {\it for any non-empty 
$C \subseteq S_\alpha$, we have 
$\sum_{i\in C} \alpha_i\ne 0$.}  

\vspace{2mm}
Note that $S_\alpha=\emptyset$ iff 
$\alpha=0$, and in this case the property 
$(\ast)$ above is satisfied vacuously.  
Therefore, we have $0\in \Omega_\ell$. 

Let $\alpha\in K^\ell$ (as above) and $B=\{u_1, u_2, ..., u_\ell\} 
\subseteq  K^n$ such that $u_i \ne u_j$ for 
all $1\le i\ne j\le \ell$. For each $f\in K[z]$, 
we set 
\begin{align} 
\alpha_{f, B}\!:&=\big(\alpha_1 f(u_1), \, 
\alpha_2 f(u_2), \, ..., \, \alpha_\ell f(u_\ell)\big)\in K^\ell.
\label{Def-af} 
\end{align}

Furthermore, we also introduce the following 
$K$-subspace of the polynomial 
algebra $K[z]$:
\begin{align} 
N_{B, \alpha}\!:=
\Big\{ f(z)\in K[z]\,\Big|\, 
\sum_{i=1}^\ell \alpha_i f(u_i)=0\Big\}.
\label{Def-NB}
\end{align}

Note that $N_{B, \alpha}$ is a co-dimension one
$K$-subspace of $K[z]$ unless $\alpha=0$, in which case 
we have $N_{B, \alpha}=K[z]$. 

\begin{lemma}\label{IdealCase} 
Let $\alpha\in K^\ell$ and $B\subseteq K^n$ be as above. 
Then $N_{B, \alpha}$ is an ideal of $K[z]$ 
iff the cardinal number $|S_\alpha|\le 1$.
\end{lemma} 

\pf The $(\Leftarrow)$ part is obvious. 
To show the $(\Rightarrow)$ part, we assume 
$|S_\alpha|\ge 2$ and derive a contradiction 
as follows. 

Without losing any generality, we assume that $\alpha_1\ne 0$. Since $|S_\alpha|\ge 2$, it is easy to see that there exists $f(z)\in N_{B, \alpha}$
such that $f(u_1)\ne 0$. Let $g(z)\in K[z]$ 
such that $g(u_1)\ne 0$ and $g(u_i)=0$ for all 
$2\le i\le n$.  Then $(fg)(u_1)\ne 0$ and 
$(fg)(u_i)=0$ for all $2\le i\le n$, whence 
$fg\not \in N_{B, \alpha}$. But this contradicts 
our assumption that  $N_{B, \alpha}$ is 
an ideal of $K[z]$.
\epfv

The lemma above determines all $\alpha\in K^\ell$ 
such that $N_{B, \alpha}$ is an ideal of $K[z]$. 
To see for which $\alpha\in K^\ell$,     
$N_{B, \alpha}$ is a Mathieu 
subspace of $K[z]$, we have the  
following result proved in 
Proposition $4.6$ 
in \cite{GIC}.

\begin{propo}\label{MS-Case} 
Let $\alpha\in K^\ell$ and $B\subseteq K^n$ be as above. 
Then $N_{B, \alpha}$ is a Mathieu subspace 
of $K[z]$ iff $\alpha \in \Omega_\ell$. 
\end{propo}

Note that the proof for the proposition above 
in \cite{GIC} is under the convenient condition 
$\alpha\in (K^\times)^\ell$. But it is 
easy to see that the proof actually 
goes through without 
this extra condition.  

The first main result of this section 
is the following proposition.

\begin{propo}\label{Nf-propo}
Let $\alpha\in K^\ell$ and $B\subseteq K^n$ be as above. 
Then we have 
\begin{align} 
\sigma(N_{B, \alpha})
&=\Big\{ f(z)\in K[z]\,\Big|\, 
|S_{\alpha_{f, B}}|\le 1 \Big\}. 
\label{Nf-propo-e1}\\
\tau(N_{B, \alpha})
&=\Big\{ f(z)\in K[z]\,\Big|\, \alpha_{f, B} 
\in \Omega_\ell\Big\}. \label{Nf-propo-e2}
\end{align}
\end{propo}
\pf Note first that for each $f(z)\in K[z]$, 
by Eqs.\,(\ref{Def-af}) and (\ref{Def-NB}) 
it is easy to see that the 
following equation holds: 
\begin{align} 
(N_{B, \alpha}:f)=N_{B, \alpha_{f, B}}.
\label{Nf-propo-pe1}
\end{align}

Then with the equation above, 
Eq.\,(\ref{Nf-propo-e1}) and 
(\ref{Nf-propo-e2}) follow immediately 
from Lemma \ref{IdealCase} 
and Proposition \ref{MS-Case}, respectively.
\epfv

\begin{rmk}\label{S6-rmk1}
Let $\alpha\in K^\ell$ such that $|S_\alpha|\ge 2$. 
Then by Eq.\,$(\ref{Nf-propo-e1})$, we see that 
the set $\sigma(N_{B, \alpha})$ of the stable elements of 
the $K$-subspace $N_{B, \alpha}\subset K[z]$, 
as contrasted to all other examples discussed in this paper, 
is actually not closed under the addition 
of the $K[z]$-module $K[z]$. 
\end{rmk}

\begin{corol}
Let $\alpha\in K^\ell$ and $B\subseteq K^n$ be as above. 
Denote by $\mathcal Z_\alpha(B)$ the ideal of $K[z]$ consisting of 
the polynomials $f(z)$ such that  
 $f(u_i)=0$ for all $i\in S_\alpha$. 
Then we have 
\begin{align} 
\mathcal Z_\alpha(B)=N_{B, \alpha}\cap 
\sigma(N_{B, \alpha})=N_{B, \alpha}\cap 
\tau(N_{B, \alpha}).
\end{align}
Consequently, $\mathcal Z_\alpha(B)$ is the 
$($unique$)$ ideal of $K[z]$ that is maximum among 
all the ideals of $K[z]$ contained 
in the $K$-subspace $N_{B, \alpha}$. 
\end{corol}
\pf First, by Eqs.\,(\ref{Def-NB}) and 
(\ref{Nf-propo-e1}), it is easy to see 
that we do have 
\begin{align} 
\mathcal Z_\alpha (B)=N_{B, \alpha}\cap 
\sigma(N_{B, \alpha}).
\end{align}
Then the corollary follows immediately from Theorem \ref{Max-submodule} 
(by viewing $K[z]$ as the $K[z]$-module in the canonical way).
\epfv

Next, we assume $K=\bC$ and $n=1$, and consider the following family of 
$\bC$-subspaces of the polynomial 
algebra $\bC[z]$ in one variable $z$. 

Let $a, b \in \bC$ with $a\ne b$.   
For each $q(z)\in \bC[z]$, we set 
\begin{align}
N_q\!:=\left\{ f\in \bC[z]\,\left|\, 
\int_a^b f(z) q(z)\, dz=0\right.\right\}. 
\label{Def-Nq}
\end{align} 

It is easy to see that $N_q$ is a co-dimension one
$\bC$-subspace of $\bC[z]$ unless $q(z)=0$, 
in which case we have $N_q=\bC[z]$. 

\begin{lemma}\label{NqCaseLemma}
For every $0\ne q(z)\in \bC[z]$, $N_q$ contains no nonzero 
ideals of $\bC[z]$. In particular, $N_q$ itself is 
not an ideal of $\bC[z]$. 
\end{lemma}
\pf Note first that by changing the variable 
$z\to (b-a)z+a$ and replacing $q(z)$ by 
$q\big( (b-a)z+a\big )$, we may assume 
that $a=0$ and $b=1$.

Now, we assume otherwise. Then there exists 
$0\ne h\in \bC[z]$ such that $(h)\subseteq N_q$, where $(h)$ 
is the ideal of $\bC[z]$ generated by $h(z)$. 
But, on the other hand, let $\bar q$ and  $\bar h$ 
be the complex conjugates of $q$ and $h$, respectively.
Then $(\bar q \bar h) (h q)$ is continuous 
and has positive values at all but finitely points 
in the interval $[0, 1]$, whence 
$\int_a^b (\bar q \bar h) (h q)\, dz > 0$. 
Consequently, we have 
$(\bar q \bar h)h \in (h)$ but $(\bar q \bar h)h \not \in N_q$, 
which is a contradiction. 
\epfv

The second main result of this section is 
the following proposition.

\begin{propo}\label{NqCase}
For each $0\ne q(z)\in \bC[z]$, we have
\begin{align} 
\sigma (N_q) &= \{0\},\label{NqCase-e1}\\
\tau (N_q)   &= N_q^c \cup \{0\}. \label{NqCase-e2}
\end{align} 
\end{propo}

In order to prove the proposition above, 
we need the following theorem proved by 
F. Pakovich \cite{Pa}. 

\begin{theo}\label{P-Thm}
Let $a\ne b \in \bC$ and $f(z), q(z)\in \bC[z]$ 
such that for all $m\ge 1$, we have 
\begin{align}
\int_a^b q(z)\, dz &\ne 0, \\
\int_a^b f^m(z) q(z)&\, dz =0.  
\end{align}
Then $f(z)=0$.
\end{theo}

\begin{corol}\label{P-Thm-Corol}
Let $0\ne q(z)\in \bC[z]$. Then  
$N_q$ is a Mathieu subspace 
of $\bC[z]$ iff $\int_a^b q(z)\, dz \ne 0$. 
\end{corol}
\pf The $(\Leftarrow)$ part follows immediately 
from Eq.\,(\ref{Def-Nq}), Theorem \ref{P-Thm} 
and Lemma \ref{radical-Lemma}. 
To show the $(\Rightarrow)$ part, 
we assume $\int_a^b q(z)\, dz = 0$, 
and derive a contradiction 
as follows.  

Note first that by Proposition \ref{Pull-Back}, it is easy to see that  
Mathieu subspaces of $\bC[z]$ are preserved by automorphisms of the polynomial algebra of $\bC[z]$. So, by applying the change of 
the variables in the proof of Lemma \ref{NqCaseLemma}, 
we may assume that $a=0$ and $b=1$.

Since $\int_0^1 q(z)\, dz = 0$, 
we have $1\in N_q$. Since $N_q$ 
by assumption is a Mathieu subspace 
of $\bC[z]$, by Lemma \ref{ONE-Lemma}  
we have $N_q=\bC[z]$.  
In particular, the complex conjugate 
$\bar q$ of $q$ lies in $N_q$, 
i.e., $\int_0^1 \bar q q \, dz = 0$.
But, this is impossible, as argued in the 
proof of Lemma \ref{NqCaseLemma}. 
Therefore, we get a contradiction. 
\epfv

\underline{\it Proof of Proposition \ref{NqCase}}\,: 
For each $h(z)\in \bC[z]$, by Eq.\,(\ref{Def-Nq}) it is easy to see 
that we have 
\begin{align} 
(N_q:h)=N_{hq}. \label{NqCase-pe1}
\end{align}

Hence by Lemma \ref{NqCaseLemma} and the equation above, 
we get Eq.\,(\ref{NqCase-e1}). 

To show  Eq.\,(\ref{NqCase-e2}), note first that 
by Lemma \ref{NqCaseLemma} and 
Theorem \ref{Max-submodule}, we have 
$\tau(N_q) \subseteq N_q^c \cup\{0\}$. 

Now, let $h(z)\in N_q^c$. 
Then we have $\int_a^b h(z) q(z)\, dz \ne 0$.
Applying Corollary \ref{P-Thm-Corol} to 
the polynomial $hq$, we see that $N_{hq}$ 
is a Mathieu subspace of $\bC[z]$. 
Then by Eq.\,(\ref{NqCase-pe1}) above, 
we get $h\in \tau(N_q)$. Hence, we have  
$N_q^c \subseteq \tau(N_q)$.  
Moreover, by Lemma \ref{EasyLemma1}, $ii)$  
we also have $0\in \tau(N_q)$, 
whence Eq.\,(\ref{NqCase-e2}) 
follows.
\epfv

{\bf Acknowledgments}\,\,The author would like to thank the referee and Michiel de Bondt for pointing out some mistakes and misprints in the earlier version of this paper.

\end{document}